\documentclass[aos,MSNbibl,nameyear,dvips]{arximspdf}
\usepackage{dcolumn}
\usepackage{graphicx}

\doi{10.1214/12-AOS1028} 
\volume{40}
\issue{4}
\pubyear{2012}
\firstpage{2043}
\lastpage{2068}

\makeatletter

\newcolumntype{d}[1]{D{.}{.}{#1}}

\renewcommand{\tilde}{\widetilde}

\newproclaim{assumption}{Condition}
\newtheorem{lemma}{Lemma}
\newtheorem{proposition}{Proposition}
\newtheorem{theorem}{Theorem}

\newproclaim{example}{Example}

\newcommand{\bu}{{\mathbf u}}
\newcommand{\bw}{{\mathbf w}}
\newcommand{\by}{{\mathbf y}}
\newcommand{\bA}{{\mathbf A}}
\newcommand{\bB}{{\mathbf B}}
\newcommand{\bT}{{\mathbf T}}
\newcommand{\bE}{{\mathbf E}}
\newcommand{\bG}{{\mathbf G}}
\newcommand{\bH}{{\mathbf H}}
\newcommand{\bI}{{\mathbf I}}
\newcommand{\bP}{{\mathbf P}}
\newcommand{\Mcal}{\mathcal{M}}
\newcommand{\Rcal}{\mathcal{R}}
\newcommand{\Gcal}{\mathcal{G}}
\newcommand{\bV}{{\mathbf V}}
\newcommand{\bX}{{\mathbf X}}
\newcommand{\bZ}{{\mathbf Z}}
\newcommand{\bzero}{{\mathbf0}}
\newcommand{\bveps}{\bolds\varepsilon}
\newcommand{\bbeta}{\bolds{\beta}}
\newcommand{\bgamma}{\bolds{\gamma}}
\newcommand{\bPsi}{\bolds\Psi}
\newcommand{\bxi}{\bolds\xi}
\newcommand{\bSig}{\bolds\Sigma}

\newcommand{\hbv}{\widehat{{\mathbf v}}}
\newcommand{\hbbeta}{\widehat\bbeta}
\newcommand{\hgamma}{\widehat\gamma}
\newcommand{\tgamma}{\widetilde\gamma}
\newcommand{\hbgamma}{\widehat\bgamma}
\newcommand{\tbgamma}{\widetilde\bgamma}
\newcommand{\hbeta}{\widehat\beta}
\newcommand{\mathcalR}{\mathcal{R}}
\newcommand{\mathcalG}{\mathcal{G}}
\newcommand{\Mfrak}{\mathfrak{M}}
\newcommand{\bMfrak}{\overline{\mathfrak{M}}}
\newcommand{\tbP}{\widetilde{\bP}}

\newcommand{\var}{\operatorname{var}}
\newcommand{\veps}{\varepsilon}
\newcommand{\Lammin}{\Lambda_{\mathrm{min}}}
\newcommand{\Lammax}{\Lambda_{\mathrm{max}}}

\makeatother

\begin{document}
\begin{frontmatter}

\title{Variable selection in linear mixed effects models}
\runtitle{Variable selection in linear mixed effects models}

\begin{aug}
\author[A]{\fnms{Yingying} \snm{Fan}\corref{}\thanksref{t1}\ead[label=e1]{fanyingy@marshall.usc.edu}}
\and
\author[B]{\fnms{Runze} \snm{Li}\thanksref{t2}\ead[label=e2]{rli@stat.psu.edu}}
\runauthor{Y. Fan and R. Li}
\affiliation{University of Southern California and Pennsylvania State University}
\address[A]{Department of Information\\
\quad and Operations Management\\
Marshall School of Business\\
University of Southern California\\
Los Angeles, California 90089\\
USA\\
\printead{e1}}
\address[B]{Department of Statistics\\
\quad and the Methodology Center\\
Pennsylvania State University\\
University Park, Pennsylvania 16802\\
USA\\
\printead{e2}} 
\end{aug}

\thankstext{t1}{Supported by NSF CAREER Award DMS-11-50318
and Grant DMS-09-06784, and 2010 Zumberge Individual Award from
USCs James H. Zumberge Faculty Research and Innovation
Fund.}

\thankstext{t2}{Supported by NIDA, NIH Grants R21 DA024260 and P50
DA10075 and in part by National Natural Science
Foundation of China Grants 11028103 and 10911120395. The content is
solely the responsibility of the authors and does not necessarily
represent the official views of the NIDA, the NNSF or the NIH.}

\received{\smonth{1} \syear{2012}}
\revised{\smonth{6} \syear{2012}}

%
\begin{abstract}
This paper is concerned with the selection and estimation of fixed and
random effects in linear mixed effects models. We propose a class of
nonconcave penalized profile likelihood methods for selecting and
estimating important fixed effects. To overcome the difficulty of
unknown covariance matrix of random effects, we propose to use a proxy
matrix in the penalized profile likelihood. We establish conditions on
the choice of the proxy matrix and show that the proposed procedure
enjoys the model selection consistency where the number of fixed
effects is allowed to grow exponentially with the sample size. We
further propose a group variable selection strategy to simultaneously
select and estimate important random effects, where the unknown
covariance matrix of random effects is replaced with a proxy matrix. We
prove that, with the proxy matrix appropriately chosen, the proposed
procedure can identify all true random effects with asymptotic
probability one, where the dimension of random effects vector is
allowed to increase exponentially with the sample size. Monte Carlo
simulation studies are conducted to examine the finite-sample
performance of the proposed procedures. We further illustrate the
proposed procedures via a real data example.
\end{abstract}

%
\begin{keyword}[class=AMS]
\kwd[Primary ]{62J05}
\kwd{62J07}
\kwd[; secondary ]{62F10}
\end{keyword}
\begin{keyword}
\kwd{Adaptive Lasso}
\kwd{linear mixed effects models}
\kwd{group variable selection}
\kwd{oracle property}
\kwd{SCAD}
\end{keyword}

\end{frontmatter}

\section{Introduction}\label{sec1}

During the last two decades, linear mixed effects models
[\citet{laird1,BLongford1}] have been widely used to model
longitudinal and repeated measurements data, and have received much
attention in the fields of agriculture, biology, economics, medicine
and sociology; see \citet{VM00} and references therein. With the
advent of modern technology, many variables can be easily collected in
a scientific study, and it is typical to include many of them in the
full model at the initial stage of modeling to reduce model
approximation error. Due to the complexity of the mixed effects models,
the inferences and interpretation of estimated models become
challenging as the dimension of fixed or random components increases.
Thus the selection of important fixed or random components becomes a
fundamental problem in the analysis of longitudinal or repeated
measurements data using mixed effects models.

Variable selection for mixed effects models has become an active
research topic in the literature.
\citet{lin0} considers testing a hypothesis on the variance component.
The testing procedures can be used to
detect whether an individual random component is significant or not.
Based on these
testing procedures, a stepwise procedure can be constructed for
selecting important random effects.
\citet{vaida1} propose the conditional AIC, an extension of the AIC
[\citet{A73}],
for mixed effects models with detailed discussion on how to define
degrees of freedom in the presence
of random effects. The conditional AIC has further been discussed in
\citet{LWZ2008}.
\citet{chen2} develop a Bayesian variable selection procedure for
selecting important random effects in the linear
mixed effects model using the Cholesky decomposition of the covariance
matrix of random effects,
and specify a prior distribution on the standard deviation of random
effects with a positive mass
at zero to achieve the sparsity of random components.
\citet{PN06} extend the generalized information criterion to select
linear mixed effects models
and study the asymptotic behavior of the proposed method for selecting
fixed effects. \citet{BG10} propose a joint variable selection method
for fixed and random effects in the linear mixed effects model using a
modified Cholesky decomposition in the setting of fixed dimensionality
for both fixed effects and random effects. \citet{IZGG2011} propose to
select fixed and random effects in a general class of mixed effects
models with fixed dimensions of both fixed and random effects using
maximum penalized likelihood method with the SCAD penalty and the
adaptive least absolute shrinkage and selection operator penalty.

In this paper, we develop a class of variable selection procedures for
both fixed effects and random effects in linear mixed effects models by
incorporating the recent advances in variable selection. We propose to
use the regularization methods to select and estimate fixed and random
effects. As advocated by \citet{fan2}, regularization methods can
avoid the stochastic error of variable selection in stepwise
procedures, and can significantly reduce computational cost compared
with the best subset selection and Bayesian procedures. Our proposal
differs from the existing ones in the literature mainly in two aspects.
First, we consider the high-dimensional setting and allow dimension of
fixed or random effects to grow exponentially with the sample size.
Second, our proposed procedures can estimate the fixed effects vector
without estimating the random effects vector and vice versa.

We first propose a class of variable selection methods for the fixed
effects using penalized
profile likelihood method. To overcome the difficulty of unknown
covariance matrix of random effects, we propose to replace it with a
suitably chosen proxy matrix. The penalized profile likelihood is
equivalent to a penalized quadratic loss
function of the fixed effects. Thus, the proposed approach can take
advantage of the recent developments
in the computation of the penalized least-squares methods [\citet
{EHJT04,ZL08}]. The optimization of the
penalized likelihood can be solved by the LARS algorithm without extra effort.
We further systematically study the sampling properties of the
resulting estimate of
fixed effects. We establish conditions on the proxy matrix and show
that the resulting estimate enjoys model selection oracle property
under such conditions. In our theoretical investigation, the number of
fixed effects is allowed to grow exponentially with the total sample
size, provided that the covariance matrix of random effects is
nonsingular. In the case of singular covariance matrix for random
effects, one can use our proposed method in Section~\ref{sec3} to first select
important random effects and then conduct variable selection for fixed
effects. In this case, the number of fixed effects needs to be smaller
than the total sample size.

Since the random effects vector is random, our main interest is in the
selection of true random effects. Observe that if a random effect
covariate is a noise variable, then the corresponding realizations of
this random effect should all be zero, and thus the random effects
vector is sparse. So we propose to first estimate the realization of
random effects vector using a group regularization method and then
identify the important ones based on the estimated random effects
vector. More specifically, under the Bayesian framework, we show that
the restricted posterior distribution of the random effects vector is
independent of the fixed effects coefficient vector. Thus, we propose a
random effect selection procedure via penalizing the restricted
posterior mode. The proposed procedure reduces the impact of error
caused by the fixed effects selection and estimation. The unknown
covariance matrix is replaced with a suitably chosen proxy matrix. In
the proposed procedure, random effects selection is carried out with
group variable selection techniques [\citet{YL06}]. The optimization of
the penalized restricted posterior mode is equivalent to
the minimization of the penalized quadratic function of random effects.
In particular, the form
of the penalized quadratic function is similar to that in the adaptive
elastic net [\citet{ZH05,ZZh09}], which allows us to minimize the
penalized quadratic function using existing algorithms. We further
study the theoretical properties of the proposed procedure and
establish conditions on the proxy matrix for ensuring the model
selection consistency of the resulting estimate. We show that, with
probability tending to\vadjust{\goodbreak} one, the proposed procedure can select all true
random effects. In our theoretical study, the dimensionality of random
effects vector is allowed to grow exponentially with the sample size as
long as the number of fixed effects is less than the total sample size.

The rest of this paper is organized as follows. Section~\ref{sec2} introduces
the penalized profile likelihood method for the estimation of fixed effects
and establishes its oracle property. We
consider the estimation of random effects and prove the model selection
consistency of the resulting
estimator in Section~\ref{sec3}. Section~\ref{sec4} provides two
simulation studies and a real data example. Some discussion is given in
Section~\ref{sec5}. All proofs are presented in Section~\ref{sec6}.

\section{Penalized profile likelihood for fixed effects}\label{sec2}

Suppose that we have a sample of $N$ subjects. For the $i$th subject,
we collect the response variable $y_{ij}$, the $d\times1$ covariate vector
${\mathbf x}_{ij}$ and $q\times1$ covariate vector ${\mathbf z}_{ij}$,
for $j = 1,\ldots, n_i$,
where $n_i$ is the number of observations on the $i$th subject. Let $n
= \sum_{i=1}^Nn_i$, $m_{n} = \max_{1\leq i\leq N}n_i$, and $\tilde
m_{n} = \min_{1\leq i\leq N}n_i$. We consider the case where $\limsup
_{n}\frac{m_{n}}{\tilde m_{n}}< \infty$, that is, the sample sizes
for $N$ subjects are balanced. For succinct presentation, we use matrix
notation and write $\by_i = (y_{i1}, y_{i2}, \ldots, y_{in_i})^T$,
$\bX_{i} = ({\mathbf x}_{i1}, {\mathbf x}_{i2},\ldots, {\mathbf
x}_{in_i})^T$ and $\bZ_i = ({\mathbf z}_{i1}, {\mathbf z}_{i2},\ldots,
{\mathbf z}_{in_i})^T$.
In linear mixed effects models, the vector of repeated measurements
$\by_i$ on the $i$th subject
is assumed to follow the linear regression model
%
\begin{equation}
\label{e001} \by_i = \bX_i \bbeta+ \bZ_i
\bgamma_i + \bveps_i,
\end{equation}
where $\bbeta$ is the $d\times1$ population-specific fixed effects
coefficient vector,
$\bgamma_i$ represents the $q\times1$ subject-specific random effects
with $\bgamma_i \sim N(\bzero, G)$,
$\bveps_i$ is the random error vector with components independent and
identically distributed as $N(0, \sigma^2)$, and $\bgamma_1,
\ldots, \bgamma_N, \bveps_1, \ldots, \bveps_N$ are independent.
Here, $G$ is the covariance matrix of random effects and may be
different from the identity matrix. So the random effects can be
correlated with each other.

Let vectors $\by$, $\bgamma$ and $\bveps$, and matrix $\bX$ be
obtained by stacking vectors $\by_i$, $\bgamma_i$ and $\bveps_i$,
and matrices $\bX_i$, respectively,
underneath each other, and let $\bZ=\operatorname{diag}\{\bZ_1,
\ldots, \bZ_N \}$ and $\mathcal{G} = \operatorname{diag}\{G,\ldots, G\}
$ be block diagonal matrices. We further standardize the design
matrix $\bX$ such that each column has norm $\sqrt{n}$. The linear
mixed effects model~(\ref{e001}) can be rewritten as
%
\begin{equation}
\label{e011} \by= \bX\bbeta+ \bZ\bgamma+ \bveps.
\end{equation}

\subsection{Selection of important fixed effects}\label{sec2.1}

In this subsection, we assume that there are no noise random effects,
and $\Gcal$ is positive definite. In the case where noise random
effects exist, one can use the method in Section~\ref{sec3} to select
the true ones. The joint density of $\by$ and $\bgamma$ is
%
\begin{eqnarray}
\label{e030}\quad f(\by, \bgamma) & = & f(\by| \bgamma) f(\bgamma) \nonumber\\
&=&(2\pi\sigma
)^{-(n+qN)/2}|\mathcal{G}|^{-1/2}
\\
&&{}\times\exp\biggl\{-\frac{1}{2\sigma^2}(\by- \bX\bbeta- \bZ
\bgamma)^T(\by- \bX\bbeta- \bZ\bgamma)-\frac{1}{2}
\bgamma^T\mathcal{G}^{-1}\bgamma\biggr\}.\nonumber
\end{eqnarray}
Given $\bbeta$, the maximum likelihood estimate (MLE) for $\bgamma$ is
$
\widehat{\gamma}(\bbeta) = \bB_z (\by-\bX\bbeta)$,
where $\bB_z = (\bZ^T \bZ+ \sigma^2\mathcal{G}^{-1}
)^{-1} \bZ^{T}$. Plugging $\widehat{\gamma}(\bbeta)$ into $f(\by,
\bgamma)$ and dropping the constant term yield the following profile
likelihood function:
%
\begin{equation}
\label{e012} L_n\bigl(\bbeta, \widehat{\bgamma}(\bbeta)\bigr)=\exp
\biggl\{-\frac
{1}{2\sigma^2}(\by- \bX\bbeta)^T\bP_z(\by-\bX
\bbeta) \biggr\},
\end{equation}
where $\bP_z= (\bI- \bZ\bB_z)^T(\bI-\bZ\bB_z) + \sigma^2\bB_z^T\mathcal
{G}^{-1}\bB_z $ with $\bI$ being the identity matrix. By
Lemma~\ref{L1} in Section~\ref{sec6}, $\bP_z$ can be rewritten
as $\bP_z = (\bI+ \sigma^{-2}\bZ\Gcal\bZ^T )^{-1}$.
To select the important $x$-variables, we propose to maximize the
following penalized profile log-likelihood function:
%
\begin{equation}
\label{e031} \log\bigl( L_n\bigl(\bbeta, \widehat{\bgamma}(\bbeta)
\bigr) \bigr) - n\sum_{j=1}^{d_n}
p_{\lambda_n} \bigl(|\beta_j|\bigr),
\end{equation}
where $p_{\lambda_n}(x)$ is a penalty function with regularization
parameter $\lambda_n\geq0$.
Here, the number of fixed effects $d_n$ may increase with sample size $n$.

Maximizing~(\ref{e031}) is equivalent to minimizing
%
\begin{equation}
\label{e010} Q_n(\bbeta) = \frac{1}{2} (\by- \bX
\bbeta)^T\bP_z(\by-\bX\bbeta)+n\sum
_{j=1}^{d_n} p_{\lambda_n} \bigl(|\beta_j|\bigr).
\end{equation}
Since $\bP_z$ depends on the unknown covariance matrix $\Gcal$ and
$\sigma^2$, we propose to use a proxy $\tbP_z = (\bI+ \bZ\Mcal\bZ
^T)^{-1}$ to replace $\bP_z$, where $\Mcal$ is a
pre-specified matrix. Denote by $\widetilde Q_n(\bbeta)$ the
corresponding objective function when $\tbP_z$ is used. We will
discuss in the next section how to choose $\Mcal$.

We note that~(\ref{e010}) does not depend on the inverse of $\Gcal$.
So although we started this section with the nonsingularity assumption
of $\Gcal$, in practice our method can be directly applied even when
noise random effects exist, as will be illustrated in simulation
studies of Section~\ref{sec4}.

Many authors have studied the selection of the penalty function to
achieve the purpose of variable selection for the linear regression
model. \citet{tibshirani3} proposes the Lasso method by the use of
$L_1$ penalty. \citet{fan2} advocate the use of nonconvex
penalties. In particular, they suggest the use of the SCAD penalty.
\citet{zou1} proposes the adaptive Lasso by using adaptive $L_1$
penalty, \citet{Zhang09} proposes the minimax concave penalty
(MCP), \citet{Liu1} propose to linearly combine $L_0$ and $L_1$
penalties and \citet{lv1} introduce a unified approach to sparse
recovery and model selection using general concave penalties. In this
paper, we use concave penalty function for variable selection.
%
\begin{assumption}\label{con2}
For each $\lambda>0$, the penalty function $p_{\lambda}(t)$ with $t
\in[0, \infty)$ is
increasing and concave with $p_{\lambda}(0) = 0$, its second order
derivative exists and is continuous and $p_{\lambda}'(0+) \in(0,
\infty)$. Further, assume that $\sup_{t>0} p_{\lambda}''(t)
\rightarrow0$ as \mbox{$\lambda\rightarrow0$}.
\end{assumption}
Condition~\ref{con2} is commonly assumed in studying regularization
methods with concave penalties. Similar conditions can be found in
\citet{fan2}, \citet{fan3} and \citet{lv1}. Although it is assumed that $p''_{\lambda
}(t)$ exists and is continuous, it can be relaxed to the case where
only $p'_{\lambda}(t)$ exists and is continuous. All theoretical
results presented in later sections can be generalized by imposing
conditions on the local concavity of $p_{\lambda}(t)$,
as in \citet{lv1}.

\subsection{Model selection consistency}\label{sec2.2}

Although the proxy matrix $\tbP_z$ may be different from the true one
$\bP_z$, solving the regularization problem~(\ref{e010}) may still
yield correct model selection results at the cost of some additional
bias. We next establish conditions on $\tbP_z$ to ensure the model
selection oracle property of the proposed method.

Let $\bbeta_0$ be the true coefficient vector. Suppose that $\bbeta_0$
is sparse, and denote $s_{1n} = \|\bbeta_0\|_0$, that is, the number of
nonzero elements in $\bbeta_0$.
Write
\[
\bbeta_0 = (\beta_{1,0}, \ldots, \beta_{d_n,0})^T
= \bigl(\bbeta_{1,0}^T, \bbeta_{2,0}^T
\bigr)^T,
\]
where $\bbeta_{1,0}$ is an $s_{1n}$-vector and $\bbeta_{2,0}$ is a
$(d_n-s_{1n})$-vector.
Without loss of generality, we assume that $\bbeta_{2,0}=\bzero$,
that is, the nonzero elements of
$\bbeta_0$ locate at the first $s_{1n}$ coordinates. With a slight
abuse of notation, we write $\bX= (\bX_1, \bX_2)$ with $\bX_1$
being a submatrix formed by the first $s_{1n}$ columns of $\bX$ and
$\bX_2$ being formed by the remaining columns. For a matrix $\bB$,
let $\Lammin(\bB)$ and $\Lammax(\bB)$ be its minimum and maximum
eigenvalues, respectively. We will need the following assumptions.
%
\begin{assumption}\label{con4}
(A) Let $a_n = \min_{1\leq j\leq s_{1n}} |\beta_{0,j}| $. It
holds that
\[
a_n n^{\tau}(\log n)^{-3/2} \rightarrow\infty
\]
with $\tau
\in(0, \frac{1}{2})$ being some positive constant, and
$\sup_{t\geq a_n/2}p''_{\lambda_n}(t) = o (n^{-1+2\tau} )$.

(B) There exists a constant $c_1>0$ such that $\Lammin
(c_1\Mcal- \sigma^{-2}\Gcal) \geq0$ and $\Lammin(c_1\sigma^{-2}(\log
n)\Gcal- \Mcal) \geq0$.

(C) The minimum and maximum eigenvalues of matrices $n^{-1}(\bX
_1^T\bX_1)$ and $n^{\theta}(\bX_1^T\bP_z\bX_1)^{-1}$ are both
bounded from below and above by $c_0$ and $c_0^{-1}$, respectively,
where $\theta\in(2\tau, 1]$ and $c_0>0$ is a constant. Further, it
holds that
%
\begin{eqnarray}
\label{e005} \biggl\| \biggl(\frac{1}{n}\bX_1^T
\tbP_z\bX_1 \biggr)^{-1} \biggr\|_\infty&\leq&
n^{-\tau}(\log n)^{3/4}/p'_{\lambda_n}(a_n/2),
\\
\label{e002}
\bigl\|\bX_2^T\tbP_z\bX_1\bigl(
\bX_1^T\tbP_z\bX_1
\bigr)^{-1}\bigr\|_\infty&<& p'_{\lambda_n}(0+)/p'_{\lambda_n}(a_n/2),
\end{eqnarray}
where $\|\cdot\|_\infty$ denotes the matrix infinity norm.
\end{assumption}

Condition~\ref{con4}(A) is on the minimum signal strength $a_n$. We
allow the minimum signal strength to decay with sample size $n$. When
concave penalties such as SCAD\vadjust{\goodbreak} [\citet{fan2}] or SICA
[\citet{lv1}] are
used, this condition can be easily satisfied with $\lambda_n$
appropriately chosen. Conditions~\ref{con4}(B) and (C) put constraints
on the proxy $\Mcal$. Condition~\ref{con4}(C) is about the design
matrices $\bX$ and $\bZ$. Inequality~(\ref{e002}) requires noise
variables and signal variables not highly correlated. The upper bound
of~(\ref{e002}) depends on the ratio $p'_{\lambda_n}(0+)/p'_{\lambda
_n}(a_n/2)$. Thus, concave penalty functions relax this condition when
compared to convex penalty functions. We will further discuss
constraints~(\ref{e005}) and~(\ref{e002}) in Lemma~\ref{L3}.

If the above conditions on the proxy matrix are satisfied, then the
bias caused by using $\tbP_z$ is small enough, and the resulting
estimate still enjoys the model selection oracle property described in
the following theorem.\vspace*{-2pt}
%
\begin{theorem}\label{T1} Assume that $\sqrt{n}\lambda_n \rightarrow
\infty$ as $n\rightarrow\infty$ and $\log d_n = o(n\lambda_n^2)$.
Then under\vspace*{1pt} Conditions~\ref{con2} and~\ref{con4}, with probability
tending to 1 as $n\rightarrow\infty$, there exists a strict local
minimizer $\hbbeta= (\hbbeta{}^T_1, \hbbeta{}^T_2)^T$ of $\tilde
Q_n(\bbeta)$ which satisfies
%
\begin{equation}
\label{e003} \|\hbbeta_1 - \bbeta_{0,1}\|_\infty<
n^{-\tau}(\log n)\quad \mbox{and}\quad \hbbeta_2 = \bzero.\vspace*{-2pt}
\end{equation}
\end{theorem}

Theorem~\ref{T1} presents the weak oracle property in the sense of
\citet{lv1} on the local minimizer of $\tilde Q(\bbeta)$. Due to
the high dimensionality and the concavity of $p_{\lambda}(\cdot)$, the
characterization of the global minimizer of $\tilde Q(\bbeta)$ is a
challenging open question. As will be shown in the simulation and real
data analysis, the concave function $\tilde Q(\bbeta)$ will be
iteratively minimized by the local linear approximation method
[\citet{ZL08}]. Following the same idea as in \citet{ZL08},
it can be shown that the resulting estimate poesses the properties in
Theorem~\ref{T1} under some conditions.\vspace*{-2pt}

\subsection{Choice of proxy matrix $\Mcal$}\label{sec2.3}

It is difficult to see from~(\ref{e005}) and~(\ref{e002}) on how
restrictive the conditions on the proxy matrix $\Mcal$ are. So we
further discuss these conditions in the lemma below. We introduce the
notation $\bT= \sigma^2\Gcal^{-1}+\bZ^T\bP_x\bZ$ and $\bE=
\sigma^2\Gcal^{-1}+\bZ^T\bZ$ with $\bP_x = \bI- \bX_1(\bX_1^T\bX
_1)^{-1}\bX_1$. Correspondingly, when the proxy matrix $\Mcal
$ is used, define $\widetilde{\bT} = \Mcal^{-1}+\bZ^T\bP_x\bZ$
and $\widetilde{\bE} = \Mcal^{-1}+\bZ^T\bZ$. We use \mbox{$\|\cdot\|_2$}
to denote the matrix $2$-norm, that is, $\|\bB\|_2 = \{\Lambda_{\mathrm{max}}(\bB
\bB^T)\}^{1/2}$ for a matrix $\bB$.\vspace*{-2pt}
%
\begin{lemma}\label{L3} Assume that $\|(\frac{1}{n}\bX_1^T\bP_z\bX
_1)^{-1}\|_\infty< n^{-\tau}\sqrt{\log n}/p'_{\lambda_n}(a_n/2)$ and
%
\begin{equation}
\label{e041} \bigl\|\bT^{-1/2}\widetilde{\bT}\bT^{-1/2} - \bI
\bigr\|_2 < \bigl( 1 + n^{\tau}s_{1n}^{1/2}p'_{\lambda_n}(a_n/2)
\bigl\|\bZ\bT^{-1}\bZ^T\bigr\|_2 \bigr)^{-1}.
\end{equation}
Then~(\ref{e005}) holds.

Similarly, assume that $\|\bX_2^T\bP_z\bX_1(\bX_1^T\bP_z\bX_1)^{-1}\|
_\infty< p'_{\lambda_n}(0+)/p'_{\lambda_n}(a_n/2)$, and
there exists a constant $c_2>0$ such that
%
\begin{eqnarray}
\label{e038}
&&
\bigl\|\bT^{-1/2}\widetilde{\bT}\bT^{-1/2} - \bI\bigr\|_2\nonumber\\
&&\qquad<
\bigl[1+n^{-1}\bigl\|\bZ\bT^{-1}\bZ^T\bigr\|_2
\\
&&\qquad\quad\hspace*{18pt}{} \times\max\bigl\{c_2n^{\theta},c_0^{-1}(
\log n)s_{1n}^{1/2}\lambda_n^{-1}p_{\lambda_n}'(a_n/2)
\bigl\|\bX_2^T\bP_z\bX_1
\bigr\|_2\bigr\} \bigr]^{-1},\nonumber
\\
\label{e045}
&&\bigl\|\bE^{-1/2}\widetilde{\bE}\bE^{-1/2}-\bI\bigr\|_2 \nonumber\\
&&\qquad<
\bigl[1+ \lambda_n^{-1}(\log n)s_{1n}^{1/2}(
\log n) p'_{\lambda
_n}(a_n/2)
\\
&&\qquad\quad\hspace*{18.5pt}{} \times\bigl\|\bZ\Gcal\bZ^T\bigr\|_2 \bigl\{\bigl\| \bigl(
\bX_1^T\bP_z \bX_1
\bigr)^{-1}\bigr\|_2\bigl\|\bX_2^T
\bP_z\bX_2 \bigr\|_2 \bigr\}^{1/2}
\bigr]^{-1},\nonumber
\end{eqnarray}
then~(\ref{e002}) holds.
\end{lemma}

Equations~(\ref{e041}),~(\ref{e038}) and~(\ref{e045}) show
conditions on the proxy matrix $\Mcal$. Note that if penalty function
used is flat outside of a neighborhood of zero, then $p'_{\lambda
_n}(a_n/\allowbreak2) \approx0$ with appropriately chosen regularization
parameter $\lambda_n$, and conditions~(\ref{e041}) and~(\ref{e045}),
respectively, reduce to
%
\begin{equation}
\label{e050} \bigl\|\bT^{-1/2}\widetilde{\bT}\bT^{-1/2} - \bI
\bigr\|_2 <1,\qquad \bigl\|\bE^{-1/2}\widetilde{\bE}\bE^{-1/2} - \bI
\bigr\|_2 < 1.
\end{equation}
Furthermore, since $\bZ$ is a block diagonal matrix, if the maximum
eigenvalue of $\bZ\bT^{-1}\bZ^T$ is of the order $o(n^{1-\theta})$,
then condition~(\ref{e038}) reduces to
%
\begin{equation}
\label{e051} \bigl\|\bT^{-1/2}\widetilde{\bT}\bT^{-1/2} - \bI
\bigr\|_2 <1.
\end{equation}
Conditions~(\ref{e050}) and~(\ref{e051}) are equivalent to assuming
that $\bT^{-1/2}\widetilde{\bT}\bT^{-1/2}$ and $\bE^{-1/2}\widetilde{\bE
}\bE^{-1/2}$ have eigenvalues bounded between 0
and 2. By linear algebra, they can further be reduced to $\|\bT
^{-1}\widetilde{\bT}\|_2 <2$ and $\|\bE^{-1}\widetilde{\bE}\|_2 <
2$. It is seen from the definitions of $\bT$, $\widetilde\bT$, $\bE
$ and $\widetilde\bE$ that if eigenvalues of $\bZ\bP_x\bZ^T$ and
$\bZ\bZ^T$ dominate those of $\sigma^2\Gcal^{-1}$ by a larger order
of magnitude, then these conditions are
not difficult to satisfy. In fact, note that both $\bZ\bP_x\bZ^T$
and $\bZ\bZ^T$ have components with magnitudes increasing with $n$,
while the components of $\sigma^{2}\Gcal^{-1}$ are independent of
$n$. Thus as long as both matrices $\bZ\bP_x\bZ^T$ and $\bZ\bZ^T$
are nonsingular, these conditions will easily be satisfied with the
choice $\Mcal= (\log n)\bI$ when $n$ is large enough.

\section{Identifying important random effects}\label{sec3}
In this section, we allow the number of random
effects $q$ to increase with sample size $n$ and write it as $q_n$ to
emphasize its dependency on $n$. We focus on the case where the number
of fixed effects $d_n$ is smaller than the total sample size $n = \sum
_{i=1}^Nn_i$. We discuss the $d_n \geq n$ case in the discussion
Section~\ref{sec5}. The major goal of this section is to select
important random effects.

\subsection{Regularized posterior mode estimate}\label{sec3.1}
The estimation of random effects is different from the estimation of
fixed effects, as the vector $\bgamma$ is random. The empirical Bayes
method has been used to estimate the random effects vector $\bgamma$
in the literature. See, for example, \citet{BBT73},
\citet{GCSR95} and \citet{VM00}.
Although the empirical Bayes method is useful in estimating random
effects in many situations, it cannot be used to select important
random effects. Moreover, the performance of an empirical Bayes
estimate largely depends on the accuracy of estimated fixed effects.
These difficulties call for a new proposal for random effects selection.

\citet{PT71} propose the error contrast method to obtain the restricted
maximum likelihood of a linear model. Following their notation, define
the $n\times(n-d)$ matrix $\bA$ by the conditions $
\bA\bA^T = \bP_x$ and $\bA^T \bA=\bI$, where $\bP_x = \bI- \bX
(\bX^T \bX)^{-1} \bX^T$. Then the vector $\bA^T\bveps$ provides a
particular set of $n-d$ linearly independent error contrasts.
Let $\bw_1 = \bA^T\by$.
The following proposition characterizes the
conditional distribution of $\bw_1$:

\begin{proposition}\label{prop1}
Given $\bgamma$, the density function of $\bw_1$ takes the form
%
\begin{equation}\quad
\label{e042} f_{\bw_1} \bigl(\bA^T\by|\bgamma\bigr) = \bigl(2
\pi\sigma^2\bigr)^{-(n-d)/2}\exp\biggl\{-\frac{1}{2\sigma^2}(\by-
\bZ\bgamma)^T\bP_x(\by-\bZ\bgamma) \biggr\}.
\end{equation}
%
\end{proposition}
The above conditional probability is independent of the fixed effects
vector $\bbeta$ and the error contrast matrix $\bA$, which allows us
to obtain a posterior mode estimate of $\bgamma$ without estimating
$\bbeta$ and calculating $\bA$.

Let $\mathfrak{M}_0\subset\{1,2,\ldots, q_n\}$ be the index set of
the true random effects. Define
\[
\bMfrak_0 =\{j\dvtx j = iq_n + k\mbox{,  for } i=0,1,2,
\ldots, N-1 \mbox{ and } k \in\Mfrak_0\}
\]
and denote by $\bMfrak_0^c = \{1, 2, \ldots, Nq_n\}\setminus\bMfrak_0$.
Then $\bMfrak_0$ is the index set of nonzero random effects
coefficients in the vector $\bgamma$, and $\bMfrak_0^c$ is the index
set of the zero ones. Let $s_{2n} = \|\Mfrak_0\|_0$ be the number of
true random effects. Then $\|\bMfrak_0 \|_0 = Ns_{2n}$.
We allow $Ns_{2n}$ to diverge with sample size $n$, which covers both
the case where the number of subjects $N$ diverges with $n$ alone and
the case where $N$ and $s_{2n}$ diverge with $n$ simultaneously.

For any $\mathcal{S}\subset\{1,\ldots, q_nN\}$, we use $\bZ_{\mathcal
{S}}$ to denote the $(q_nN)\times|\mathcal{S}|$ submatrix
of $\bZ$ formed by columns with indices in $\mathcal{S}$, and
$\bgamma_{\mathcal{S}}$ to denote the subvector of $\bgamma$ formed
by components with indices in $\mathcal{S}$. Then $\bgamma_{\bMfrak
_0} \sim N(\bzero, \Gcal_{\bMfrak_0})$ with $\Gcal_{\bMfrak_0}$ a
submatrix formed by entries of $\Gcal$ with row and column indices in
$\bMfrak_0$. In view of~(\ref{e042}), the restricted posterior
density of $\bgamma_{\bMfrak_0}$ can be derived as
\begin{eqnarray*}
&&f_{\bw_1}\bigl(\bgamma_{\bMfrak_0} |\bA^T\by\bigr)\propto
f_{\bw_1}\bigl(\bA^T\by|\bgamma_{\bMfrak_0}\bigr)f(
\bgamma_{\bMfrak_0})
\\
&&\qquad\propto\exp\biggl\{-\frac{1}{2\sigma^2}(\by-\bZ_{\bMfrak
_0}
\bgamma_{\bMfrak_0})^T \bP_x(\by-\bZ_{\bMfrak_0}
\bgamma_{\bMfrak_0})-\frac{1}{2}\bgamma_{\bMfrak_0}^T
\mathcal{G}_{\bMfrak
_0}^{-1}\bgamma_{\bMfrak_0} \biggr\}.
\end{eqnarray*}
Therefore, the restricted posterior mode estimate of $\bgamma_{\bMfrak
_0}$ is the solution to the following minimization problem:
%
\begin{equation}
\label{e007} \min_{\bgamma} \bigl\{(\by-\bZ_{\bMfrak_0}
\bgamma_{\bMfrak_0})^T \bP_x(\by-\bZ_{\bMfrak_0}
\bgamma_{\bMfrak_0}) + \sigma^2\bgamma_{\bMfrak_0}^T
\mathcal{G}_{\bMfrak_0}^{-1}\bgamma_{\bMfrak_0} \bigr\}.
\end{equation}

In practice, since the true random effects $\bMfrak_0$ are unknown,
the formulation~(\ref{e007}) does not help us estimate $\bgamma$. To
overcome this difficulty, note that $\bZ_{\bMfrak_0}\bgamma_{\bMfrak
_0}=\bZ\bgamma$ and $\bgamma_{\bMfrak_0}^T\mathcal{G}_{\bMfrak
_0}^{-1}\bgamma_{\bMfrak_0} = \bgamma^T\mathcal{G}^{+}\bgamma$
with $\Gcal^{+}$ the Moore--Penrose generalized inverse of~$\Gcal$.
Thus the objective function in~(\ref{e007}) is rewritten as
\[
(\by-\bZ\bgamma)^T \bP_x(\by-\bZ\bgamma) +
\sigma^2\bgamma^T\mathcal{G}^{+}\bgamma,
\]
which no longer depends on the unknown $\bMfrak_0$.
Observe that if the $k$th random effect is a noise one, then the
corresponding standard deviation is 0, and the coefficients
$\gamma_{ik}$ for all subjects $i=1, \ldots, N$ should equal to 0.
This leads us to consider group\vspace*{1pt} variable selection strategy to identify
true random effects. Define $\gamma_{\cdot k}=(\sum_{i=1}^N\gamma
_{ik}^2)^{1/2}$,
$k=1,\ldots, q_n$, and consider the following regularization problem:
%
\begin{equation}
\label{e008} \frac{1}{2}(\by-\bZ\bgamma)^T \bP_x(
\by-\bZ\bgamma)+\frac
{1}{2}\sigma^2 \bgamma^T
\mathcal{G}^{+}\bgamma+ n\sum_{k=1}^{q_n}
p_{\lambda_n}(\gamma_{\cdot k}),
\end{equation}
where $p_{\lambda_n}(\cdot)$ is the penalty function with
regularization parameter $\lambda_n\geq0$. The penalty function here may
be different from the one in Section~\ref{sec2}. However, to ease the
presentation, we use the same notation.

There are several advantages to estimating the random effects vector
$\bgamma$ using the above proposed method~(\ref{e008}). First, this
method does not require knowing or estimating the fixed effects vector
$\bbeta$, so it is easy to implement, and the estimation error of
$\bbeta$ has no impact on the estimation of $\bgamma$. In addition,
by using the group variable selection technique, the true random
effects can be simultaneously selected and estimated.

In practice, the covariance matrix $\Gcal$ and the variance $\sigma^2$
are both unknown. Thus, we replace $\sigma^{-2}\mathcal{G}$ with
$\Mcal$, where $\Mcal=\operatorname{diag}\{M,\ldots, M\}$ with $M$
a proxy of~$G$, yielding the following regularization problem:
%
\begin{equation}\label{e006}\quad
\tilde{Q}_n^*(\bgamma)=\frac{1}{2}(\by-\bZ
\bgamma)^T \bP_x(\by-\bZ\bgamma)+\frac{1}{2}
\bgamma^T\Mcal^{-1}\bgamma+ n \sum
_{k=1}^{q_n} p_{\lambda_n}(\gamma_{\cdot k}).
\end{equation}
It is interesting to observe that the form of regularization in (\ref
{e006}) includes the elastic net [\citet{ZH05}] and the adaptive
elastic net [\citet{ZZh09}] as special cases. Furthermore, the
optimization algorithm for adaptive elastic net can be modified for
minimizing~(\ref{e006}).

\subsection{Asymptotic properties}\label{sec3.2}

Minimizing~(\ref{e006}) yields an estimate of $\bgamma$, denoted by
$\widehat\bgamma$.
In this subsection, we study the asymptotic property of $\widehat
{\bgamma}$.
Because $\bgamma$ is random rather than a deterministic parameter vector,
the existing formulation for the asymptotic analysis of a
regularization problem
is inapplicable to our setting. Thus, asymptotic analysis
of $\widehat{\bgamma}$ is challenging.

Let $\bT= \bZ^T\bP_x\bZ+ \sigma^2\Gcal^{+}$ and $\widetilde\bT=
\bZ^T\bP_x\bZ+ \Mcal^{-1}$.
Denote by $\bT_{11} = \bZ^T_{\bMfrak_0}\bP_x\*\bZ_{\bMfrak_0} +
\sigma^2(\Gcal_{\bMfrak_0})^{-1}$, $\bT_{22}= \bZ^T_{\bMfrak
_0^c}\bP_x\bZ_{\bMfrak_0^c}$ and $\bT_{12} = \bZ^T_{\bMfrak_0}\bP_x\bZ
_{\bMfrak_0^c}$. Similarly, we\vadjust{\goodbreak} can define submatrices
$\widetilde\bT_{11}$, $\widetilde\bT_{22}$ and $\widetilde\bT_{12}$ by
replacing $\sigma^{-2}\Gcal$ with $\Mcal$. Then it is easy
to see that $\widetilde\bT_{12}=\bT_{12}$. Notice that if the oracle
information of set ${\bMfrak_0}$ is available and $\Gcal$, and
$\sigma^2$ are known, then the Bayes estimate of the true random
effects coefficient vector $\bgamma_{\bMfrak_0}$ has the form $\bT
_{11}^{-1}\bZ_{\bMfrak_0}^T\bP_x\by$. Define $\bgamma^* =
((\bgamma_1^*)^T,\ldots,(\bgamma_N^*)^T)^T$ with $\bgamma_j^* =
(\gamma_{j1}^*,\ldots,\gamma_{jq_n}^*)^T$ for $j = 1, \ldots, N$ as
the oracle-assisted Bayes estimate of the random effects vector. Then
${\bgamma^*}_{\bMfrak_0^c} = \bzero$ and ${\bgamma^*}_{\bMfrak_0}
= \bT_{11}^{-1}\bZ_{\bMfrak_0}^T\bP_x\by$. Correspondingly, define
$\tbgamma^*$ as the oracle Bayses estimate with proxy matrix, that is,
$\tbgamma^*_{\bMfrak_0^c} = \bzero$ and
%
\begin{equation}
\label{e016} {\tbgamma^*}_{\bMfrak_0} = \widetilde\bT_{11}^{-1}
\bZ_{\bMfrak
_0}^T\bP_x\by.
\end{equation}
For $k=1,\ldots, q_n$, let $\gamma_{\cdot k}^* = \{\sum_{j=1}^N(\gamma
_{jk}^*)^2\}^{1/2}$. Throughout we condition on the event
%
\begin{equation}
\label{e026} \Omega^* = \Bigl\{\min_{k\in\mathfrak{M}_0} \gamma_{\cdot
k}^*\geq
\sqrt{N}b_0^*\Bigr\}
\end{equation}
with $b_0^* \in(0, \min_{\in\Mfrak_0}\sigma_k)$ and $\sigma_k^2 =
\var(\gamma_{jk})$. The above event $\Omega^*$ is to ensure that the
oracle-assisted estimator $\gamma_{\cdot k}^*/\sqrt{N}$ of $\sigma_k$
is not too negatively biased. 

\begin{assumption}\label{con3}
(A) The maximum eigenvalues satisfy $\Lammax(\bZ_iG\bZ_i^T)\leq
c_3s_{2n}$ for all $i=1,\ldots, N$ and the minimum and
maximum eigenvalues of
$m_n^{-1}\bZ_{\bMfrak_0}^T\bP_x \bZ_{\bMfrak_0}$ and $G_{\Mfrak
_0}$ are bounded from below and above by $c_3$ and $c_3^{-1}$,
respectively, with $m_n = \max_{1\leq i\leq N}n_i$, where $c_3$ is a
positive constant. Further, assume that for some $\delta\in(0, \frac{1}{2})$,
%
\begin{eqnarray}
\label{e017}
\bigl\|\widetilde\bT_{11}^{-1}\bigr\|_{\infty}&\leq&
\frac{\sqrt{N}n^{-1-\delta}}{p'_{\lambda_n} (\sqrt{N}b_0^*/2 )},
\\
\label{e044}
\max_{j \in\Mfrak_0^c}\bigl\|\tilde\bZ_j^T\bP_x
\bZ_{\bMfrak
_0}\widetilde\bT_{11}^{-1}\bigr\|_2 &<&
\frac{p'_{\lambda
_n}(0+)}{p'_{\lambda_n}(\sqrt{N}b_0^*/2)},
\end{eqnarray}
where $\tilde\bZ_j$ is the submatrix formed by the $N$ columns of
$\bZ$ corresponding to the $j$th random effect.

(B) It holds that $\sup_{\{t \geq\sqrt{N}b_0^*/2\}}
p''_{\lambda_n}(t) = o(N^{-1})$.

(C) The proxy matrix satisfies $\Lammin(\Mcal- \sigma^{-2}\Gcal)
\geq0$.
\end{assumption}

Condition~\ref{con3}(A) is about the design matrices $\bX$, $\bZ$
and covariance matrix~$\Gcal$. Since $\bZ_{\bMfrak_0}$ is a block
diagonal matrix and $\limsup\frac{\max_{i}n_i}{\min_{i}n_i} <
\infty$, the components of $\bZ_{\bMfrak_0}^T\bP_x\bZ_{\bMfrak
_0}$ have magnitude of the order $m_n=O(n/N)$. Thus, it is not very
restrictive to assume that the minimum and maximum eigenvalues of $\bZ
_{\bMfrak_0}^T\bP_x\bZ_{\bMfrak_0}$ are both of the order $m_n$.
Condition~(\ref{e044}) puts an upper bound on the correlation between
noise covariates and true covariates. The upper bound of~(\ref{e044})
depends on the penalty function. Note that for concave penalty we have
$p'_{\lambda_n}(0+)/p'_{\lambda_n}(\sqrt{N}b_0^*/2) >1$, whereas\vadjust{\goodbreak} for
$L_1$ penalty \mbox{$p'_{\lambda_n}(0+)/p'_{\lambda_n}(\sqrt{N}b_0^*/2) =
1$}. Thus,\vspace*{1pt} concave penalty relaxes~(\ref{e044}) when compared with the
$L_1$ penalty. Condition~\ref{con3}(B) is satisfied by many commonly
used penalties with appropriately chosen $\lambda_n$, for example,
$L_1$ penalty, SCAD penalty and SICA penalty with small $a$. Condition
\ref{con3}(C) is a restriction on the proxy matrix $\Mcal$, which
will be further discussed in the next subsection.

Let $\bgamma= (\bgamma_{1}^T, \ldots, \bgamma_{N}^T)^T$ with
$\bgamma_{j}= (\gamma_{j1},\ldots, \gamma_{jq_n})^T$ being an
arbitrary $(Nq_n)$-vector. Define $\gamma_{\cdot k} = (\sum
_{j=1}^{N}\gamma_{jk}^2 )^{1/2}$ for each $k=1, \ldots, q_n$. Let
%
\begin{equation}
\label{e023} \mathfrak{M} (\bgamma) = \bigl\{k\in\{1,\ldots, q_n\}\dvtx
\gamma_{\cdot k} \neq0 \bigr\}.
\end{equation}
Theorem~\ref{T2} below shows that there exists a local minimizer of
$\widetilde{Q}_n^*(\bgamma)$ defined in~(\ref{e006}) whose support
is the same as the true one $\bMfrak_0$, and that this local minimizer
is close to the oracle estimator $\tbgamma^*$.
%
\begin{theorem}\label{T2}
Assume that Conditions~\ref{con2} and~\ref{con3} hold,
$b_0^*n^{\delta}/\sqrt{N} \rightarrow\infty$, $\log(Nq_n
)=o (n^2\lambda_n^2/(Ns_{2n}m_{n}) )$, and $n^2\lambda
_n^2/(Nm_{n}s_{2n}) \rightarrow\infty$ as $n\rightarrow\infty$.
Then, with probability tending to 1, there exists a strict local
minimizer $\hbgamma\in\mathbf{R}^{Nq_n}$ of $\widetilde
{Q}_n^*(\bgamma)$ such that
\[
\mathfrak{M}(\hbgamma) = \mathfrak{M}_0 \quad\mbox{and}\quad
\max_{k\in
\Mfrak_0} \Biggl\{\frac{1}{N}\sum_{j=1}^{N}
\bigl(\hgamma_{jk}-\tgamma_{jk}^*\bigr)^2 \Biggr
\}^{1/2} \leq n^{-\delta},
\]
where $\delta$ is defined in~(\ref{e017}).
\end{theorem}

Using a similar argument to that for Theorem~\ref{T1}, we can obtain
that the dimensionality $Nq_n$ is also allowed to grow exponentially
with sample size $n$ under some growth conditions and with
appropriately chosen $\lambda_n$. In fact, note that if the sample
sizes $n_1 =\cdots=n_N \equiv m_{n}/N$, then the growth condition in
Theorem~\ref{T2} becomes $\log(Nq_n) = o(ns_{2n}^{-1}\lambda_n^2)$.
Since the lowest signal level in this case is $\sqrt{N}b_0^*$, if
$b_0^*$ is a constant, a reasonable choice of tuning parameter would be
of the order $\sqrt{N}n^{-\kappa}$ with some $\kappa\in(0, \frac
{1}{2})$. For $s_{2n}=O(n^{\nu})$ with $\nu\in[0, \frac{1}{2})$ and
$Nn^{1-2\kappa-\nu} \rightarrow\infty$, we obtain that $Nq_n$ can
grow with rate $\exp(Nn^{1-2\kappa-\nu})$.

\subsection{Choice of proxy matrix $\Mcal$}\label{sec3.3}

Similarly as for the fixed effects selection and estimation, we discuss
(\ref{e017}) and~(\ref{e044}) in the following lemma.
%
\begin{lemma} \label{L2}
Assume that $\|\bT_{11}^{-1}\|_{\infty}< \frac{\sqrt{N}n^{-1-\delta
}}{p'_{\lambda_n}(\sqrt{N} b_0^*/2)}[1- \frac{1}{\sqrt{\log n}}]$ and
%
\begin{equation}\label{eqproxy}
\bigl\|\bT_{11}^{-1}\widetilde{\bT}_{11}-\bI
\bigr\|_2 \leq\bigl[1+\sqrt{ s_{2n}\log n}n^{1+\delta}p'_{\lambda_n}
\bigl(\sqrt{N} b_0^*/2\bigr)\bigl\|\bT_{11}^{-1}
\bigr\|_2 \bigr]^{-1}.
\end{equation}
Then~(\ref{e017}) holds.

Assume that $\max_{j \in\Mfrak_0^c}\|\tilde\bZ_j^T\bP_x\bZ_{\bMfrak
_0}\bT_{11}^{-1}\|_2 < \frac{p'_{\lambda_n}(0+)}{2
p'_{\lambda_n}(\sqrt{N}b_0^*/2)}$ with $\tilde\bZ_j$ defined in
(\ref{e044}) and
%
\begin{equation}\label{eqproxy1}
\bigl\|\bT_{11}\widetilde\bT_{11}^{-1}-\bI
\bigr\|_2 \leq1.
\end{equation}
Then~(\ref{e044}) holds.
\end{lemma}

Conditions~(\ref{eqproxy}) and~(\ref{eqproxy1}) put restrictions on
the proxy matrix $\Mcal$. Similarly to the discussions after Lemma
\ref{L3}, if $p'_{\lambda_n}(\sqrt{N}b_0^*/2) \approx0$, then these
conditions become $\|\bT_{11}\widetilde\bT_{11}^{-1} - \bI\|_2 <1$.
If $\bZ_{\bMfrak_0}^T\bP_x\bZ_{\bMfrak_0}$ dominates $\sigma^{2}\Gcal
_{\bMfrak_0}^{-1}$ by a larger magnitude, then conditions
(\ref{eqproxy}) and~(\ref{eqproxy1}) are not restrictive, and
choosing $\Mcal= (\log n)\bI$ should make these conditions as well as
Condition~\ref{con3}(C) satisfied for large enough~$n$.

We remark that using the proxy matrix $\Mcal=(\log n)\bI$ is equivalent
to ignoring correlations among random effects. The idea of using
diagonal matrix as a proxy of covariance matrix has been proposed in
other settings of high-dimensional statistical inference. For instance,
the naive Bayes rule (or independence rule), which replaces the full
covariance matrix in Fisher's discriminant analysis with a diagonal
matrix, has been demonstrated to be advantageous for high-dimensional
classifications both theoretically [\citet{Bickel04,FF08}] and
empirically [\citet{Dudoit2002}]. The intuition is that although
ignoring correlations gives only a biased estimate of covariance
matrix, it avoids the errors caused by estimating a large amount of
parameters in covariance matrix in high dimensions. Since the
accumulated estimation error can be much larger than the bias, using
diagonal proxy matrix indeed produces better results.

\section{Simulation and application}
\label{sec4} In this section, we investigate the finite-sample
performance of the proposed procedures by simulation studies and a real
data analysis. Throughout, the SCAD penalty with $a=3.7$
[\citet{fan2}] is used. For each simulation study, we randomly simulate 200
data sets. Tuning parameter selection plays an important role in
regularization methods. For fixed effect selection, both AIC- and
BIC-selectors [\citet{ZhangLiTsai10}] are used to select the
regularization parameter $\lambda_n$ in~(\ref{e010}). Our simulation
results clearly indicate that the BIC-selector performs better than the
AIC-selector for both the SCAD and the LASSO penalties. This is
consistent with the theoretical analysis in \citet{WangLiTsai07}.
To save space, we report the results with the BIC-selector. Furthermore
the BIC-selector is used for fixed effect selection throughout this
section. For random effect selection, both AIC- and BIC-selectors are
also used to select the regularization parameter $\lambda_n$ in
(\ref{e006}). Our simulation results imply that the BIC-selector
outperforms the AIC-selector for the LASSO penalty, while the SCAD with
AIC-selector performs better than the SCAD with BIC-selector. As a
result, we use AIC-selector for the SCAD and BIC-selector for the LASSO
for random effect selection throughout this section.
\begin{example}\label{Example1}
We compare our method with some existing ones in
the literature under the same model setting as that in
\citet{BG10}, where a joint variable selection method for fixed
and random effects in linear mixed effects models is proposed. The
underlying true model takes the following form with $q = 4$ random
effects and $d = 9$ fixed effects:
%
\begin{equation}
\label{e070} \qquad
y_{ij} = b_{i1} + \beta_1x_{ij1}
+ \beta_2x_{ij2} + b_{i2}z_{ij1} +
b_{i3}z_{ij2} + \veps_{ij}, \veps_{ij}
\sim_{\mathrm{i.i.d.}} N(0,1 ),
\end{equation}
where the true parameter vector $\bbeta_0 = (1,1,0,\ldots,0)^T$, the
true covariance matrix for random effects
\[
\bG= \pmatrix{
9 & 4.8 & 0.6
\cr
4.8 & 4 & 1
\cr
0.6 & 1 & 1
}
\]
and the covariates $x_{ijk}$ for $k=1,\ldots, 9$ and $z_{ijl}$ for
$l=1,2,3$ are generated independently from a uniform distribution over
the interval $[-2,2]$. So there are three true random effects and two
true fixed effects. Following \citet{BG10}, we consider two different
sample sizes $N=30$ subjects and $n_i = 5$ observations per subject,
and $N = 60$ and $n_i = 10$. Under this model setting, \citet{BG10}
compared their method with various methods in the literature, and
simulations therein demonstrate that their method outperforms the
competing ones. So we will only compare our methods with the one in
\citet{BG10}.

In implementation, the proxy matrix is chosen as $\Mcal= (\log n) \bI
$. We then estimate the fixed effects vector $\bbeta$ by minimizing
$\widetilde{Q}_n(\bbeta)$, and the random effects vector
$\bgamma$ by minimizing~(\ref{e006}). To understand the effects of
using proxy matrix $\Mcal$ on the estimated
random effects and fixed effects, we compare our estimates with the
ones obtained by solving regularization problems~(\ref{e010}) and
(\ref{e008}) with the true value $\sigma^{-2}\Gcal$.

\begin{table}
\tablewidth=235pt
\caption{Fixed and random effects selection in Example \protect\ref{Example1} when~$d=9$
and~$q=4$}\label{tab1}
\begin{tabular*}{\tablewidth}{@{\extracolsep{\fill}}lcd{3.1}d{3.1}@{}}
\hline
\textbf{Setting} & \textbf{Method} & \multicolumn{1}{c}{\textbf{\%CF}}
& \multicolumn{1}{c@{}}{\textbf{\%CR}} \\
\hline
$N=30$ & Lasso-P & 51 & 19.5 \\
$n_i=5$ & SCAD-P & 90 & 86 \\
& SCAD-T & 93.5 & 99 \\
& M-ALASSO & 73 & 79 \\
[6pt]
$N=60$ & Lasso-P & 52 &50.5 \\
$n_i=10$ & SCAD-P & 100 & 100 \\
& SCAD-T & 100 & 100\\
& M-ALASSO & 83 & 89 \\
\hline
\end{tabular*}
\end{table}

Table~\ref{tab1} summarizes the results by using our method with the
proxy matrix $\Mcal$ and SCAD penalty (SCAD-P),
our method with proxy matrix $\Mcal$ and Lasso penalty (Lasso-P), our
method with true $\sigma^{-2}\Gcal$ and SCAD penalty
(SCAD-T). When SCAD penalty is used, the local linear approximation
(LLA) method proposed by \citet{ZL08} is employed
to solve these regularization problems. The rows ``M-ALASSO'' in Table
\ref{tab1} correspond to the joint estimation
method by \citet{BG10} using BIC to select the tuning parameter. As
demonstrated in \citet{BG10}, the BIC-selector
outperforms the AIC-selector for M-ALASSO. We compare these methods by
calculating the percentage of times the correct fixed
effects are selected (\%CF), and the percentage of times the correct
random effects are selected (\%CR).
Since these two measures were also used in \citet{BG10}, for simplicity
and fairness of comparison, the results for
M-ALASSO in Table~\ref{tab1} are copied from \citet{BG10}.

It is seen from Table~\ref{tab1} that SCAD-P greatly outperforms
Lasso-P and M-ALASSO. We also see that when the true covariance matrix
$\sigma^{-2}\Gcal$ is used, \mbox{SCAD-T} has almost perfect variable
selection results. Using the proxy matrix makes the results slightly
inferior, but the difference vanishes for larger sample size $N=60, n_i=10$.
\end{example}
\begin{example}\label{Example2}
In this example, we consider the case where
the design matrices for fixed and random effects overlap. The sample
size is fixed at $n_i=8$ and $N=30$, and the numbers for fixed and
random effects are chosen to be $d=100$ and $q=10$, respectively. To
generate the fixed effects design matrix, we first independently generate
$\tilde{{\mathbf x}}_{ij}$ from $N_d(\bzero, \bSig)$, where $\bSig
=(\sigma_{st})$ with
$\sigma_{st}=\rho^{|s-t|}$ and $\rho\in(-1,1)$. Then for the $j$th
observation of the $i$th subject, we set $x_{ijk}=I(\tilde
{x}_{ijk}>0)$ for covariates $k=1$ and $d$, and
set $x_{ijk}=\tilde{x}_{ijk}$ for all other values of $k$. Thus 2 out
of $d$ covariates are discrete ones and the rest are continuous ones.
Moreover, all covariates are correlated with each other. The covariates
for random effects are the same as the corresponding ones for fixed
effects, that is, for the $j$th observation of the $i$th subject, we
set $z_{ijk} = x_{ijk}$ for $ k=1,\ldots, q=10$. Then the random
effect covariates form a subset of fixed effect covariates.

The first six elements of fixed effects vector $\bbeta_0$ are
$(2,0,1.5, 0,0,1)^T$, and the remaining elements are all zero. The
random effects vector $\bgamma$ is generated in the same way as in
Example~\ref{Example1}. So the first covariate is discrete and has both nonzero
fixed and random effect. We consider different values of correlation
level $\rho$, as shown in Table~\ref{tab7}. We choose $\Mcal= (\log
n)\bI$.

\begin{table}[b]
\def\arraystretch{0.9}
\caption{Fixed and random effects selection and estimation in Example
\protect\ref{Example2} when $n_i=8$, $N=30$, $d=100$, $q=10$ and design matrices for fixed
and random effects overlap}\label{tab7}\label{tabaddlabel}%
\begin{tabular*}{\tablewidth}{@{\extracolsep{\fill}}lcd{2.2}d{2.2}cc
d{2.2}ccc@{}}
\hline
& & \multicolumn{4}{c}{\textbf{Random effects}} &
\multicolumn{4}{c@{}}{\textbf{Fixed effects}} \\[-4pt]
& & \multicolumn{4}{c}{\hrulefill} &
\multicolumn{4}{c@{}}{\hrulefill} \\
&&\multicolumn{1}{c}{\textbf{FNR}}
& \multicolumn{1}{c}{\textbf{FPR}} & \multicolumn{1}{c}{\textbf{MRL2}}
& \multicolumn{1}{c}{\textbf{MRL1}} & \multicolumn{1}{c}{\textbf{FNR}}
& \multicolumn{1}{c}{\textbf{FPR}} & \multicolumn{1}{c}{\textbf{MRL2}}
& \multicolumn{1}{c@{}}{\textbf{MRL1}} \\
\textbf{Setting} & \textbf{Method} &  \multicolumn{1}{c}{\textbf{(\%)}} & \multicolumn{1}{c}{\textbf{(\%)}} & &
& \multicolumn{1}{c}{\textbf{(\%)}}
& \multicolumn{1}{c}{\textbf{(\%)}} & &\\
\hline
$\rho=0.3$ & Lasso-P & 11.83 & 9.50& 0.532&0.619& 62.67 & 0.41 & 0.841
& 0.758\\
& SCAD-P & 0.50 & 1.07 & 0.298 & 0.348 & 0.83 & 0.03 & 0.142 & 0.109 \\
& SCAD-T & 3.83 & 0.00 & 0.522 & 0.141 & 0.33 & 0.02 & 0.102 & 0.082 \\
[6pt]
$\rho=-0.3$ & Lasso-P & 23.67 & 7.64 & 0.524 & 0.580& 59.17 & 0.41 &
0.802 & 0.745 \\
& SCAD-P & 1.83 & 0.71 & 0.308 & 0.352 & 0.67 & 0.05 & 0.141 & 0.109 \\
& SCAD-T & 3.17 & 0.00 & 0.546 & 0.141 & 0.17 & 0.02 & 0.095 & 0.078 \\
[6pt]
$\rho=0.5$ & Lasso-P & 9.83 & 10.07 & 0.548 & 0.631 & 60.33 & 0.48 &
0.844 & 0.751 \\
& SCAD-P & 1.67 & 0.50 & 0.303 & 0.346 & 0.17 & 0.05 & 0.138 & 0.110 \\
& SCAD-T & 5.00 & 0.00 & 0.532 & 0.149 & 0.50 & 0.02 & 0.113 & 0.091\\
\hline
\end{tabular*}
\end{table}

Since the dimension of random effects vector $\bgamma$ is much larger
than the total sample size, as suggested at the beginning of Section
\ref{sec2.1}, we start with the random effects selection by first choosing a
relatively small tuning parameter $\lambda$ and use our method in
Section~\ref{sec3} to select important random effects. Then with the
selected random effects, we apply our method in Section~\ref{sec2} to
select fixed effects. To improve the selection results for random
effects, we further use our method in Section~\ref{sec3} with the
newly selected fixed effects to reselect random effects. This iterative
procedure is applied to both Lasso-P and SCAD-P methods. For SCAD-T,
since the true $\sigma^{-2}\Gcal$ is used, it is unnecessary to use
the iterative procedure, and we apply our methods only once for both
fixed and random effects selection and estimation.

We evaluate each estimate by calculating the relative $L_2$ estimation loss
\[
\mathrm{RL}2(\hbbeta) = \|\hbbeta- \bbeta_0\|_2/ \|
\bbeta_0\|_2,
\]
where $\hbbeta$ is an estimate of the fixed effects vector $\bbeta_0$.
Similarly, the relative $L_1$ estimation error of $\hbbeta$,
denoted by $\mathrm{RL}1(\hbbeta)$, can be calculated by replacing the
$L_2$-norm with the $L_1$-norm. For the random effects estimation, we
define $\mathrm{RL}2(\hbgamma)$ and $\mathrm{RL}1(\hbgamma)$ in a similar way by
replacing $\bbeta_0$ with the true $\bgamma$ in each simulation. We
calculate the mean values of $\mathrm{RL}2$ and $\mathrm{RL}1$ in the simulations and
denote them by MRL2 and MRL1 in Table~\ref{tab7}.
In addition to mean relative losses, we also calculate the percentages
of missed true covaritates (FNR), as well as the percentages of falsely
selected noise covariates (FPR), to evaluate the performance of
proposed methods.

\begin{table}[b]
\tabcolsep=0pt
\caption{The estimated coefficients of fixed and random effects in
Example \protect\ref{Example3}}\label{tab6}
\begin{tabular*}{\tablewidth}{@{\extracolsep{\fill}}ld{2.3}d{1.3}d{1.3}d{1.3}
d{2.3}d{1.3}d{1.3}cc@{}}
\hline
& \multicolumn{1}{c}{\textbf{Intercept}} & \multicolumn{1}{c}{$\bolds{b_1(t)}$}
& \multicolumn{1}{c}{$\bolds{b_2(t)}$} & \multicolumn{1}{c}{$\bolds{b_3(t)}$}
& \multicolumn{1}{c}{$\bolds{b_4(t)}$} & \multicolumn{1}{c}{$\bolds{b_5(t)}$} &
\multicolumn{1}{c}{$\bolds{x_1}$} & \multicolumn{1}{c}{$\bolds{x_2}$}
& \multicolumn{1}{c@{}}{$\bolds{x_3}$} \\
\hline
Fixed & 29.28 & 9.56& 5.75 & 0 & -8.32& 0 & 4.95& 0 & 0\\
Random & 0 & 0& 0& 0 &0& 0 &0 &0& 0 \\
\hline
& \multicolumn{1}{c}{$\bolds{b_1(t)x_3}$}& \multicolumn{1}{c}{$\bolds{b_2(t)x_3}$}
& \multicolumn{1}{c}{$\bolds{b_3(t)x_3}$}
& \multicolumn{1}{c}{$\bolds{b_4(t)x_3}$} & \multicolumn{1}{c}{$\bolds{b_5(t)x_3}$}&
\multicolumn{1}{c}{$\bolds{x_1x_2}$}& \multicolumn{1}{c}{$\bolds{x_1x_3}$}& \multicolumn{1}{c}{$\bolds{x_2x_3}$}
& \\
\hline
Fixed & 0& 0& 0 &0& 0 &0 &0 & 0&\\
Random & 0.163 & 0.153 & 0.057 & 0.043& 0.059 & 0 & 0.028 &0.055&\\
\hline
\end{tabular*}
\end{table}

From Table~\ref{tab7} we see that SCAD-T has almost perfect variable
selection results for fixed effects, while SCAD-P has highly comparable
performance, for all\vadjust{\goodbreak} three values of correlation level $\rho$. Both
methods greatly outperform the \mbox{Lasso-P} method. For the random effects
selection, both SCAD-P and SCAD-T perform very well with SCAD-T having
slightly larger false negative rates. We remark that the superior
performance of SCAD-P is partially because of the iterative procedure.
In these high-dimensional settings, directly applying our random
effects selection method in Section~\ref{sec3} produces slightly
inferior results to the ones for SCAD-T in Table~\ref{tab7}, but
iterating once improves the results. We also see that as the
correlation level increases, the performance of all methods become
worse, but the SCAD-P is still comparable to SCAD-T, and both perform
very well in all settings.
\end{example}
\begin{example}\label{Example3}
We illustrate our new procedures through an empirical
analysis of a subset of data collected in the Multi-center AIDs Cohort Study.
Details of the study design, method and medical implications have been
given by \citet{kaslow87}.
This data set comprises the human immunodeficiency virus (HIV) status
of 284 homosexual men
who were infected with HIV during the follow-up period between 1984 and 1991.
All patients are scheduled to take measurements semiannually. However,
due to the missing of scheduled visits and the random occurrence of HIV
infections,
there are an unequal number of measurements and different measurement
times for each patients.
The total number of observations is 1765.

Of interest is to investigate the relation between the mean CD4
percentage after the infection
($y$) and predictors smoking status ($x_1$, 1 for smoker and 0 for nonsmoker),
age at infection ($x_2$), and pre-HIV infection CD4 percentage (Pre-CD4
for short, $x_3$). To account
for the effect of time, we use a five-dimensional cubic spline
${\mathbf b}(t)=(b_1(t), b_2(t), \ldots, b_5(t))^T$.
We take into account the two-way interactions ${\mathbf
b}(t_{ij})x_{i3}$, $x_{i1}x_{i2}$, $x_{i1}x_{i3}$ and
$x_{i2}x_{i3}$. These eight interactions together with variables
${\mathbf b}(t_{ij})$, $x_{i1}$, $x_{i2}$
and $x_{i3}$ give us 16 variables in total. We use these 16 variables
together with an intercept to
fit a mixed effects model with dimensions for fixed and random effects
$d=q =17$.
The estimation results are listed in Table~\ref{tab6} with rows
``Fixed'' showing the estimated $\beta_j$'s for fixed effects, and rows
``Random'' showing the estimates\vadjust{\goodbreak} $\gamma_{\cdot k}/\sqrt{N}$.
The standard error for the null model is 11.45, and it reduces to 3.76
for the selected model. From Table~\ref{tab6}, it can be seen that the
baseline has time-variant fixed effect and Pre-CD4 has time-variant
random effect. Smoking has fixed effect while age and Pre-CD4 have no
fixed effects. The interactions smoking${}\times{}$Pre-CD4 and age${}\times{}
$Pre-CD4 have random effects with smallest standard deviations among
selected random effects. The boxplots of the selected random effects
are shown in Figure~\ref{fig1}.

\begin{figure}

\includegraphics{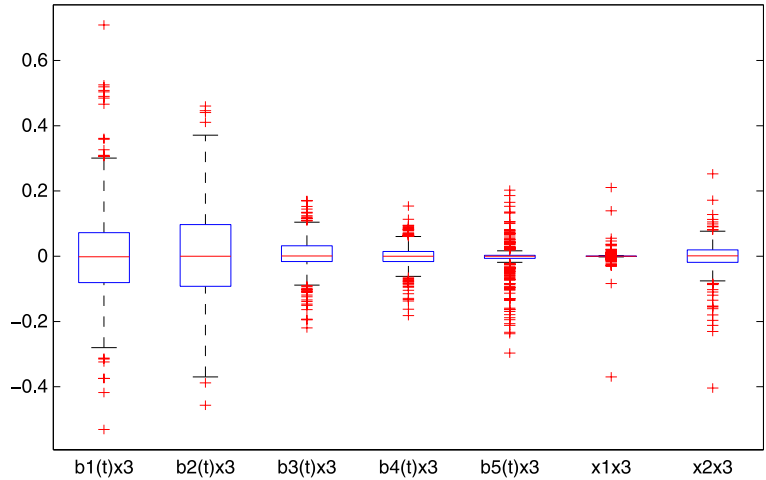}

\caption{Boxplots of selected random effects. From left to right:
$b_i(t)x_3$, $i=1,2,\ldots,5$, $x_1x_3$, $x_2x_3$, where $x_1$ is the
smoking status, $x_2$ is the age at infection, $x_3$ is Pre-CD4 level
and $b_i(t)$'s are cubic spline basis functions of time.}
\label{fig1}
\end{figure}

Our results have close connections with the ones in \citet
{Huangetal2002} and \citet{QL06}, where the former used bootstrap
approach to test the significance of variables and the later proposed
hypothesis test based on penalized spline and quadratic inference
function approaches, for varying-coefficient models. Both papers
revealed significant evidence for time-varying baseline, which is
consistent with our discovery that basis functions $b_j(t)$'s have
nonzero fixed effect coefficients. At 5\% level,
\citet{Huangetal2002} failed to reject the hypothesis of constant
Pre-CD4 effect ($p$-value 0.059), while Qu and Li's (\citeyear{QL06})
test was weakly significant with $p$-value 0.045. Our results show that
Pre-CD4 has constant fixed effect and time-varying random effect, which
may provide an explanation on the small difference of $p$-values in
\citet{Huangetal2002} and \citet{QL06}.

To further access the significance of selected fixed effects, we refit
the linear mixed effects model with selected fixed and random effects
using the Matlab function ``nlmefit.'' Based on the $t$-statistics from
the refitted model, the intercept, the baseline functions $b_1(t)$ and
$b_2(t)$ are all highly significant with $t$-statistics much larger
than 7, while the $t$-statistics for $b_4(t)$ and $x_1$ (smoking) are
$-$1.026 and 2.216, respectively. This indicates that $b_4(t)$ is
insignificant, and smoking is only weakly significant at 5\%
significance level. This result is different from those in \citet
{Huangetal2002} and \citet{QL06}, where neither paper found significant
evidence for smoking. A possible explanation is that by taking into
account random effects and variable selection, our method has better
discovery power.\vspace*{-2pt}
\end{example}

\section{Discussion}\label{sec5}

We have discussed the selection and estimation of fixed effects in
Section~\ref{sec2}, providing that the random
effects vector has nonsingular covariance matrix, while we have
discussed the selection of
random effects in Section~\ref{sec3}, providing that the dimension of fixed
effects vector is smaller than the sample size.
We have also illustrated our methods with numerical studies. In
practical implementation, the
dimensions of the random effects vector and fixed effects vector can be both
much larger than the total sample size. In such case, we suggest an
iterative way to select and estimate
the fixed and random effects. Specifically, we can first start with the
fixed effects selection using the
penalized least squares by ignoring all random effects to reduce the
number of fixed effects to below
sample size. Then in the second step, with the selected fixed effects,
we can apply our new method in
Section~\ref{sec3} to select important random effects. Third, with the
selected random effects
from the second step, we can use our method in Section~\ref{sec2} to
further select important fixed effects.
We can also iterate the second and third steps several times to improve
the model selection and estimation results.\vspace*{-2pt}

\section{Proofs}\label{sec6}

Lemma~\ref{L1} is proved in the supplemental article \citet{FL12}.\vspace*{-2pt}
%
\begin{lemma}\label{L1} It holds that
\[
\bP_z=(I - \bZ\bB_z)^T
\Rcal^{-1}(I-\bZ\bB_z) + \bB_z^T
\mathcal{G}^{-1}\bB_z=\bigl(\mathcalR+ \bZ\mathcalG
\bZ^T\bigr)^{-1}.\vspace*{-2pt}
\]
\end{lemma}


\subsection{\texorpdfstring{Proof of Theorem \protect\ref{T1}}{Proof of Theorem 1}}\label{sec6.1}

Let $\mathcal{N}_0 = \{\bbeta= (\bbeta_1^T, \bbeta_2^T)^T\dvtx \|\bbeta_1 -
\bbeta_{0,1}\|_{\infty}\leq\break n^{-\tau}(\log n),
\bbeta_2 = \bzero\in\mathbf{R}^{d_n-s_{1n}} \}$. We are going to
show that under Conditions~\ref{con2} and~\ref{con4}, there exists a
strict local minimizer $\hbbeta\in\mathcal{N}_0$ of $\tilde
Q_n(\bbeta)$ with asymptotic probability one.

For a vector $\bbeta= (\beta_1, \ldots, \beta_p)^T$, let $\bar
{p}'_{\lambda_n}(\bbeta)$ be a vector of the same length whose $j$th
component is $p_{\lambda_n}'(|\beta_j|)\operatorname{sgn}(\beta_j)$,
$j=1,\ldots, d_n$. By \citet{lv1}, the sufficient conditions for
$\hbbeta= (\hbbeta{}^T_1, \bzero^T)^T \in\mathbf{R}^{d_n}$ with
$\hbbeta_1 \in\mathbf{R}^{s_{1n}}$ being a strict local minimizer of
$\tilde Q_n(\bbeta)$ are
%
\begin{eqnarray}
\label{e058}
&\displaystyle
-\bX_1^T\tbP_z(\by- \bX_1
\hbbeta_1) + n\bar{p}'_{\lambda
_n}(\hbbeta_1) = 0,&
\\
\label{e059}
&\displaystyle \|{\mathbf v}_2\|_\infty< np_{\lambda_n}'(0+),&
\\
\label{e060}
&\displaystyle \Lammin\bigl(\bX_1^T\tbP_z
\bX_1\bigr)> -n p_{\lambda_n}''\bigl(|
\hbeta_j|\bigr),\qquad j=1,\ldots, s_{1n},&
\end{eqnarray}
where ${\mathbf v}_2 = \bX_2^T\tbP_z(\by-\bX_1\hbbeta_1)$. So we
only need to show that with probability tending to 1, there exists a
$\hbbeta\in\mathcal{N}_0$ satisfying
conditions~(\ref{e058})--(\ref{e060}).\vadjust{\goodbreak}

We first consider~(\ref{e058}). Since $\by= \bX_1\bbeta_{0,1} + \bZ
\bgamma+ \bveps$, equation~(\ref{e058}) can be rewritten as
%
\begin{equation}\quad
\label{e064} \hbbeta_1 - \bbeta_{0,1} = \bigl(
\bX_1^T\tbP_z\bX_1
\bigr)^{-1}\bX_1^T\tbP_z(\bZ
\bgamma+ \bveps) - n\bigl(\bX_1^T\tbP_z
\bX_1\bigr)^{-1}\bar{p}'_{\lambda_n}(
\hbbeta_1).
\end{equation}
Define a vector-valued continuous function
\[
g(\bbeta_1) = \bbeta_1 - \bbeta_{0,1} -
\bigl(\bX_1^T\tbP_z\bX_1
\bigr)^{-1} \bX_1^T \tbP_z (\bZ
\bgamma+ \bveps) + n\bigl(\bX_1^T\tbP_z
\bX_1\bigr)^{-1} \bar{p}'_{\lambda_n}(
\bbeta_1)
\]
with\vspace*{1pt} $\bbeta_1 \in\mathbf{R}^{s_{1n}}$. It suffices to show that
with probability tending to 1, there exists $\hbbeta=(\hbbeta{}^T_1,
\hbbeta{}^T_2)^T \in\mathcal{N}_0$ such that $g(\hbbeta_1) = 0$. To
this end, first note that
\fontsize{10pt}{\baselineskip}\selectfont
{\[
\bigl(\bX_1^T\tbP_z \bX_1
\bigr)^{-1}\bX_1^T\tbP_z(\bZ
\bgamma+ \bveps) \sim N\bigl(0, \bigl(\bX_1^T
\tbP_z \bX_1\bigr)^{-1}\bX_1^T
\tbP_z\bP_z^{-1}\tbP_z
\bX_1\bigl(\bX_1^T\tbP_z
\bX_1\bigr)^{-1}\bigr).
\]}
\normalsize
By Condition~\ref{con4}(B), the matrix
$c_1\tbP_z\! -\! \tbP_z\bP_z^{-1}\tbP_z \!=\! \tbP_z\bZ(c_1\Mcal\!-\! \sigma
^{-2}\Gcal)\bZ^T\tbP_z \!\geq0$, where $A\geq0$ means the matrix $A$
is positive semi-definite. Therefore,
%
\begin{equation}
\label{e071} \quad\bV\equiv\bigl(\bX_1^T\tbP_z
\bX_1\bigr)^{-1}\bX_1^T
\tbP_z\bP_z^{-1}\tbP_z
\bX_1 \bigl(\bX_1^T\tbP_z
\bX_1\bigr)^{-1} \leq c_1\bigl(
\bX_1^T\tbP_z \bX_1
\bigr)^{-1}.
\end{equation}
Thus, the $j$th diagonal component of matrix $\bV$ in~(\ref{e071}) is
bounded from above by the $j$th diagonal component of $c_1(\bX_1^T\tbP_z
\bX_1)^{-1}$. Further note that by Condition~\ref{con4}(B), $\tbP
_z^{-1} - c_1(\log n)\bP_z^{-1} \leq\bZ(\Mcal- c_1\frac{(\log
n)}{\sigma^2}\Gcal)\bZ^T \leq0$. Recall that by linear algebra, if
two positive definite matrices $A$ and $B$ satisfy $A \geq B$, then it
follows from the Woodbury formula that $A^{-1} \leq B^{-1}$. Thus,
$(c_1\log n)\tbP_z \geq\bP_z$ and $(\bX_1^T\tbP_z\bX_1)^{-1} \leq
(c_1\log n)(\bX_1^T\bP_z\bX_1)^{-1}$.
So by Condition~\ref{con4}(C), the diagonal components of $\bV$ in
(\ref{e071}) are bounded from above by $O(n^{-\theta}(\log n))$. This
indicates that the variance of each component of the normal random
vector $(\bX_1^T\tbP_z\bX_1)^{-1}\bX_1^T\tbP_z(\bZ\bgamma+
\bveps) $ is bounded from above by $O(n^{-\theta}(\log n))$. Hence,
by Condition~\ref{con4}(C),
%
\begin{eqnarray}
\label{e061}
\bigl\|\bigl(\bX_1^T
\tbP_z\bX_1\bigr)^{-1}\bX_1
\tbP_z(\bZ\bgamma+ \bveps)\bigr\|_\infty
&=& O_p\bigl(n^{-\theta/2}\sqrt{(\log n) (\log
s_{1n})}\bigr)\nonumber\\[-8pt]\\[-8pt]
&=& o_p \bigl(n^{-\tau}(\log n)
\bigr).\nonumber
\end{eqnarray}
Next, by Condition~\ref{con4}(A), for any $\bbeta=
(\beta_1,\ldots, \beta_{d_n})^T \in\mathcal{N}_0$ and large enough~$n$,
we can
obtain that
%
\begin{equation}
\label{e063} |\beta_j|\geq|\beta_{0,j}| - |
\beta_{0,j} - \beta_{j}|\geq a_n/2,\qquad j=1,
\ldots, s_{1n}.
\end{equation}
Since $p'_{\lambda_n}(x)$ is a decreasing function in $(0,\infty)$,
we have $\|\bar{p}'_{\lambda_n}(\bbeta_1)\|_\infty\leq p'_{\lambda
_n}(a_n/2)$. This together with Condition~\ref{con4}(C) ensures that
%
\begin{eqnarray}
\label{e062} \bigl\|\bigl(\bX_1^T \tbP_z
\bX_1\bigr)^{-1}\bar{p}'_{\lambda_n}(
\bbeta_1)\bigr\|_\infty&\leq&\bigl\|\bigl(\bX_1^T
\tbP_z\bX_1\bigr)^{-1}\bigr\|_\infty\bigl\|\bar
{p}'_{\lambda_n}(\bbeta_1)\bigr\|_\infty\nonumber\\[-8pt]\\[-8pt]
&\leq& o
\bigl(n^{-\tau-1}(\log n) \bigr).\nonumber
\end{eqnarray}
Combining~(\ref{e061}) and~(\ref{e062}) ensures that with probability
tending to 1, if $n$ is large enough,
\[
\bigl\|\bigl(\bX_1^T\tbP_z\bX_1
\bigr)^{-1}\bX_1\tbP_z(\bZ\bgamma+ \bveps) + n
\bigl(\bX_1^T\tbP_z\bX_1
\bigr)^{-1}\bar{p}'_{\lambda}(\bbeta_1)
\bigr\|_\infty< n^{-\tau}(\log n).
\]
Applying\vspace*{1pt} Miranda's existence theorem [\citet{v1989}] to the function
$g(\bbeta_1)$ ensures that there exists a vector $\hbbeta_1\in
\mathbf{R}^{s_{1n}}$ satisfying $\|\hbbeta_1 - \bbeta_{0,1}\|_\infty
< n^{-\tau}\log n$ such that $g(\hbbeta_1) = 0$.

Now we prove that the solution to~(\ref{e058}) satisfies (\ref
{e059}). Plugging $\by= \bX_1\bbeta_{0,1} + \bZ\bgamma+ \bveps$
into ${\mathbf v}$ in~(\ref{e059}) and by~(\ref{e064}), we obtain that
\[
{\mathbf v}_2 = \bX_2^T\tbP_z
\bX_1(\bbeta_{0,1}-\hbbeta_1)+
\bX_2^T\tbP_z(\bZ\bgamma+ \bveps) = {\mathbf
v}_{2,1} + {\mathbf v}_{2,2},
\]
where ${\mathbf v}_{2,1} = [-\bX_2^T\tbP_z\bX_1(\bX_1^T\tbP_z\bX
_1)^{-1}\bX_1^T\tbP_z + \bX_2^T\tbP_z](\bZ\bgamma+ \bveps)$ and
${\mathbf v}_{2,2} = \bX_2^T\tbP_z\*\bX_1(\bX_1^T\tbP_z\bX_1)^{-1}\bar
{p}_{\lambda_n}(\hbbeta_1)$.
Since $(\bZ\bgamma+ \bveps)\sim N(0, \bP_z^{-1})$, it is easy to
see that
${\mathbf v}_{2,1}$ has normal distribution with mean 0 and variance
\[
\bX_2^T \bigl( \bI- \tbP_z
\bX_1\bigl(\bX_1^T\tbP_z
\bX_1\bigr)^{-1}\bX_1^T \bigr)
\tbP_z\bP_z^{-1} \tbP_z \bigl(
\bI- \bX_1\bigl(\bX_1^T\tbP_z
\bX_1\bigr)^{-1}\bX_1^T
\tbP_z \bigr)\bX_2.
\]
Since $\bP_z^{-1}\leq c_1\tbP_z^{-1}$, $\bI- \tbP_z^{1/2}\bX_1(\bX
_1^T\tbP_z\bX_1)^{-1}\bX_1^T\tbP_z^{1/2}$ is a projection matrix,
and $\tbP_z$ has eigenvalues less than 1, it follows that for the unit
vector ${\mathbf e}_k$,
\begin{eqnarray*}
{\mathbf e}_k^T\var({\mathbf v}_{2,1}){
\mathbf e}_k &\leq& c_1{\mathbf e}_k^T
\bX_2^T \bigl( \tbP_z - \tbP_z
\bX_1\bigl(\bX_1^T\tbP_z
\bX_1\bigr)^{-1}\bX_1^T
\tbP_z \bigr)\bX_2{\mathbf e}_k
\\
&\leq& c_1{\mathbf e}_k^T
\bX_2^T\tbP_z\bX_2{\mathbf
e}_k\leq{\mathbf e}_k^T
\bX_2^T\bX_2{\mathbf e}_k =
c_1n,
\end{eqnarray*}
where in the the last step, each column of $\bX$ is standardized to
have $L_2$-norm $\sqrt{n}$.
Thus the diagonal elements of the covariance matrix of ${\mathbf
v}_{1,2}$ are bounded from above by $c_1n$. Therefore, for some large
enough constant $C>0$,
\begin{eqnarray*}
P \bigl(\|{\mathbf v}_{2,1}\|_\infty\geq\sqrt{2Cn \log
d_n} \bigr)&\leq& (d_n-s_{1n})P \bigl(\bigl|N(0,
c_1n)\bigr|\geq\sqrt{2Cn\log d_n} \bigr)
\\
&=& (d_n - s_{1n})\exp\bigl(-c_1^{-1}C
\log d_n\bigr) \rightarrow0.
\end{eqnarray*}
Thus, it follows from the assumption $\log d_n = o(n\lambda_n^2)$ that
\[
\|{\mathbf v}_{2,1}\|_\infty= O_p\bigl(\sqrt{n\log d_n}\bigr) = o_p
\bigl(np'_{\lambda_n}(0+) \bigr).
\]
Moreover, by Conditions~\ref{con4}(B) and (C),
\[
\|{\mathbf v}_{2,2}\|_\infty\leq n\bigl\|\bX_2^T
\tbP_z\bX_1\bigl(\bX_1^T
\tbP_z\bX_1\bigr)^{-1}\bigr\|_\infty
p'_{\lambda_n}(a_n/2) < np'_{\lambda_n}(0+).
\]
Therefore inequality~(\ref{e059}) holds with probability tending to 1
as $n\rightarrow\infty$.

Finally we prove that $\hbbeta\in\mathcal{N}_0$ satisfying (\ref
{e058}) and~(\ref{e059}) also makes~(\ref{e060}) hold with
probability tending to 1. By~(\ref{e063}) and Condition~\ref{con4}(A),
\[
0\leq-n p_{\lambda_n}''\bigl(|\hbeta_j|\bigr)
\leq-n \sup_{t \geq
a_n/2}p_{\lambda_n}''(t) = o
\bigl(n^{2\tau} \bigr).
\]
On the other hand, by Condition~\ref{con4}(C), $\Lammin(\bX_1^T\tbP_z\bX
_1) \geq c_0n^{\theta}$. Since $\theta> 2\tau$, inequality
(\ref{e062}) holds with probability tending to 1 as $n\rightarrow
\infty$.

Combing the above results, we have shown that with probability tending
to 1 as $n \rightarrow\infty$, there exists $\hbbeta\in\mathcal
{N}_0$ which is a strict local minimizer of $\tilde Q_n(\bbeta)$. This
completes the proof.\vadjust{\goodbreak}



\subsection{\texorpdfstring{Proof of Theorem \protect\ref{T2}}{Proof of Theorem 2}}\label{sec6.2}

Let $\bgamma= (\bgamma_1^T, \ldots, \bgamma_N^T)^T \in\mathbf
{R}^{q_nN}$ with $\bgamma_{j}^T= (\gamma_{j1},\allowbreak\ldots, \gamma_{jq_n})$
be a $\mathbf{R}^{Nq_n}$-vector satisfying $\mathfrak
{M}(\bgamma) = \mathfrak{M}_0$. Define $\bu(\bgamma) =
(\bu_1^T,\ldots,\allowbreak
\bu_N^T)^T \in\mathbf{R}^{Nq_n}$ with $\bu_j
=(u_{j1}, \ldots, u_{jq_n})^T $, where for $j=1,\ldots, N$,
%
\begin{equation}
\label{e009} \lambda_n u_{jk}=p'_{\lambda_n}(
\gamma_{\cdot k})\gamma_{jk}/\gamma_{\cdot k} \qquad\mbox{if } k\in
\mathfrak{M}(\bgamma)
\end{equation}
and\vspace*{1pt} $\lambda_n u_{jk} = 0$ if $k\notin\mathfrak{M}(\bgamma)$. Here,
$\gamma_{\cdot k} = \{\sum_{j=1}^N \gamma_{jk}^2\}^{1/2}$. Let
$\tbgamma^*$ be the oracle-assisted estimate defined in~(\ref{e016}).
By \citet{lv1}, the sufficient conditions for $\bgamma$ with $\bgamma
_{\bMfrak_0^c} = \bzero$ being a strict local minimizer of (\ref
{e006}) are
%
\begin{eqnarray}
\label{e013}
{\tbgamma^*}_{\bMfrak_0} - \bgamma_{\bMfrak_0} &=& n \lambda_n
\widetilde\bT_{11}^{-1}\bu(\bgamma_{\bMfrak_0}),
\\
\label{e014}
\Biggl(\sum_{j=1}^N
w_{jk}^2 \Biggr)^{1/2} &<& np'_{\lambda
_n}(0+),\qquad
k\in\mathfrak{M}_0^c,
\\
\label{e015}
\Lammin(\widetilde\bT_{11}) &>& n\Lammax\Biggl(-\frac{\partial
^2}{\partial\bgamma_{\bMfrak_0}^2}
\Biggl(\sum_{j=1}^{q_n} p_{\lambda_n}(
\gamma_{\cdot k}) \Biggr) \Biggr),
\end{eqnarray}
where $\bw(\bgamma) = (\bw_1^T,\ldots, \bw_N^T)^T \in\mathbf
{R}^{Nq_n}$ with $\bw_j =(w_{j1}, \ldots, w_{jq_n})^T $, and
%
\begin{equation}
\label{e018} \bw(\bgamma) = \bZ^T\bP_x(\by- \bZ\bgamma)
- \Mcal^{-1}\bgamma.
\end{equation}
We will show that, under Conditions~\ref{con2} and~\ref{con3},
conditions~(\ref{e013})--(\ref{e015}) above are satisfied with
probability tending to 1 in a small neighborhood of $\tbgamma^*$.

In general, it is not always guaranteed that~(\ref{e013}) has a
solution. We first show that under Condition~\ref{con3}, there exists
a vector $\hbgamma^*$ with $\mathfrak{M}(\hbgamma^*) = \mathfrak
{M}_0$ such that $\hbgamma^*_{\bMfrak_0}$ makes~(\ref{e013}) hold.
To this end, we constrain the objective function $\widetilde
{Q}_n^*(\bgamma)$ defined in~(\ref{e006}) on the
$(Ns_{n2})$-dimensional subspace $\mathcal{B} = \{\bgamma\in\mathbf
{R}^{q_nN}\dvtx \bgamma_{\bMfrak_0^c} = 0 \}$ of $\mathbf{R}^{q_nN}$.
Next define
\[
\mathcal{N}_1 = \Biggl\{\bgamma\in\mathcal{B}\dvtx \max_{k\in\Mfrak_0}
\Biggl\{\sum_{j=1}^{N}\bigl(
\gamma_{jk}-\tgamma_{jk}^*\bigr)^2 \Biggr
\}^{1/2} \leq\sqrt{N}n^{-\delta} \Biggr\}.
\]
For any $\tbgamma= (\tgamma_{11}, \ldots, \tgamma_{1q_n}, \ldots,
\tgamma_{N1}, \ldots,\gamma_{Nq_n})^T\in\mathcal{N}_1$ and $k\in
\Mfrak_0$, we have
%
\begin{eqnarray}\label{e028}
\bigl\|\tbgamma- \tbgamma^*\bigr\|_\infty&\leq&\max_{k\in\Mfrak_0} \Biggl\{ \sum
_{j=1}^{N}\bigl(\gamma_{jk}-
\tgamma_{jk}^*\bigr)^2 \Biggr\}^{1/2} \leq
\sqrt{N}n^{-\delta} \quad\mbox{and }
\nonumber\\
\tgamma^*_{\cdot k} &=& \Biggl\{\sum_{j=1}^N
\bigl(\tgamma^*_{jk}\bigr)^2 \Biggr\}^{1/2} \leq
\Biggl\{\sum_{j=1}^N\bigl(
\tgamma^*_{jk}-\tgamma_{jk}\bigr)^2 \Biggr
\}^{1/2} + \Biggl\{\sum_{j=1}^N(
\tgamma_{jk})^2 \Biggr\}^{1/2} \\
&\leq&\sqrt
{N}n^{-\delta} + \tgamma_{\cdot k}.\nonumber
\end{eqnarray}
Note that by Condition~\ref{con3}(C), we have $\widetilde\bT_{11}^{-1}
\geq\bT_{11}^{-1}$. Thus it can be derived using linear
algebra and the definitions of $\tgamma_{\cdot k}^*$ and $\gamma_{\cdot
k}^*$ that $\tgamma_{\cdot k}^* \geq\gamma_{\cdot k}^*$.\vadjust{\goodbreak}
Since we condition on the event $\Omega^*$ in~(\ref{e026}), it is
seen that for large enough $n$,
%
\begin{equation}
\label{e056} \tgamma_{\cdot k } \geq\tgamma^*_{\cdot k}-
\sqrt{N}n^{-\delta}\geq\gamma^*_{\cdot k}-\sqrt{N}n^{-\delta} >
\sqrt{N}b_0^*/2
\end{equation}
for $k\in\Mfrak_0$ and $\tbgamma\in\mathcal{N}_1$. Thus, in view
of the definition of $\bu(\bgamma)$ in~(\ref{e009}), for $k \in
\mathfrak{M}_0$, we have
\[
\bigl\|\lambda_n\bu(\tbgamma_{\bMfrak_0})\bigr\|_\infty\leq
\max_{k\in\mathfrak{M}_0}p'_{\lambda_n}(\widetilde{
\gamma}_{\cdot
k})\leq p'_{\lambda_n}\bigl(
\sqrt{N}b_0^*/2\bigr),
\]
where in the last step, $p'_{\lambda_n}(t)$ is decreasing in $t\in
(0,\infty)$ due to the concavity of $p_{\lambda_n}(t)$. This together
with~(\ref{e017}) in Condition~\ref{con3} ensures
%
\begin{equation}
\label{e029} \bigl\|n\lambda_n\widetilde\bT_{11}^{-1}
\bu(\tbgamma_{\bMfrak_0})\bigr\|_{\infty} \leq n\bigl\|\widetilde
\bT_{11}^{-1}\bigr\|_\infty p'_{\lambda
_n}
\bigl(\sqrt{N}b_0^*/2\bigr)\leq\sqrt{N}n^{-\delta}.
\end{equation}
Now define the vector-valued continuous function
$\bPsi(\bxi)\!=\! \bxi\!-\! {\tbgamma^*}_{\bMfrak_0}\! -\! n\lambda_n\widetilde\bT
_{11}^{-1}\bu(\bxi)
$, with $\bxi$ a $\mathbf{R}^{Ns_{2n}}$-vector.
Combining~(\ref{e028}) and~(\ref{e029}) and applying Miranda's
existence theorem [\citet{v1989}] to the function $\bPsi(\bxi)$, we
conclude that there exists $\hbgamma^* \in\mathcal{N}_1$ such that
$\hbgamma^*_{\bMfrak_0}$ is a solution to equation~(\ref{e013}).

We next show that $\hbgamma^*$ defined above indeed satisfies (\ref
{e015}). Note that for any vector ${\mathbf x}\neq\bzero$,
%
\begin{equation}
\label{e021} \frac{\partial^2}{\partial{\mathbf x}^2}p_{\lambda_n}\bigl(\|
{\mathbf x}\|_2\bigr) =
p''_{\lambda_n}\bigl(\|{\mathbf x}\|_2\bigr)
\frac{{\mathbf x}{\mathbf
x}^T}{\|{\mathbf x}\|_2}+ p'_{\lambda_n}\bigl(\|{\mathbf x}\|_2\bigr)
\biggl(\frac
{1}{\|{\mathbf x}\|_2}-\frac{{\mathbf x}{\mathbf x}^T}{\|{\mathbf x}\|
_2^3}\biggr).
\end{equation}
Since $-p'_{\lambda_n}(t)\leq0$ and $-p''_{\lambda_n}(t)\geq0$ for
$t\in(0,\infty)$, we have
\begin{eqnarray*}
\Lammax\biggl(-\frac{\partial^2}{\partial{\mathbf x}^2} p_{\lambda
_n}\bigl(\| {\mathbf x}
\|_2\bigr) \biggr) &\leq&-p''_{\lambda_n}\bigl(\|{
\mathbf x}\|_2\bigr) + \frac{p'_{\lambda_n}(\|
{\mathbf x}\|_2)}{\|{\mathbf x}\|_2} - \frac{p'_{\lambda_n}(\|
{\mathbf x}\|_2)}{\|{\mathbf x}\|_2} \\
&=&
-p''_{\lambda_n}\bigl(\|{\mathbf x}\|_2\bigr).
\end{eqnarray*}
Since $\hbgamma^*_{\bMfrak_0} \in\mathcal{N}_1$, by~(\ref{e056})
we have $\hgamma_{\cdot k}^*>\sqrt{N}b_0^*/2$ for $k \in\Mfrak_0$.
It follows from the above inequality and Condition~\ref{con3}(B) that
with probability tending to 1, the maximum eigenvalue of the matrix
$-\frac{\partial^2}{\partial\bgamma_{\bMfrak_0}^2} (\sum_{j=1}^{q_n}
p_{\lambda_n}(\hgamma_{\cdot k}^*) )$ is less than
\[
\max_{j\in\Mfrak_0} \bigl(-p''_{\lambda_n}\bigl(
\hgamma_{\cdot
j}^*\bigr) \bigr) =o\bigl(N^{-1}\bigr) =
o(m_n/n).
\]
Further, by Condition~\ref{con3}(A), $\frac{1}{n}\Lammin(\widetilde
\bT_{11}) = \frac{1}{n} \Lammin(\bZ_{\bMfrak_0}^T\bP_x\bZ_{\bMfrak_0})
\geq c_3\frac{m_{n}}{n}
$. Thus
the maximum eigenvalue of the matrix $-\frac{\partial^2}{\partial
\bgamma_{\bMfrak_0}^2} (\sum_{j=1}^{q_n}
p_{\lambda_n}(\hgamma_{\cdot k}^*) )$ is less than $n^{-1}\Lammin
(\widetilde\bT_{11})$ with asymptotic probability 1, and (\ref
{e015}) holds for $\hbgamma^*$.

It remains to show that $\hbgamma^*$ satisfies~(\ref{e014}). Let
$\hbv= \hbgamma^* - \tbgamma^*$. Since\vspace*{2pt} $\hbgamma^* $ is a solution
to~(\ref{e013}), we have
$\hbv= n \lambda_n \widetilde\bT_{11}^{-1}\bu(\hbgamma_{\bMfrak_0})$.
In view of~(\ref{e018}), we have
%
\begin{eqnarray}
\label{e025}
\bw\bigl(\hbgamma^*_{\bMfrak_0^c}\bigr) &=& \bigl(
\bZ_{\bMfrak
_0^c}^T - \widetilde\bT_{12}^T
\widetilde\bT_{11}^{-1}\bZ_{\bMfrak
_0}^T \bigr)
\bP_x\by+ \widetilde\bT_{12}^T
\hbv_{\bMfrak_0}
\nonumber\\
&=& \bigl(\bZ_{\bMfrak_0^c}^T - \widetilde\bT_{12}^T
\widetilde\bT_{11}^{-1}\bZ_{\bMfrak_0}^T \bigr)
\bP_x(\bZ\bgamma+ \bveps)+ \widetilde\bT_{12}^T
\hbv_{\bMfrak_0}\\
&\equiv&\tilde\bw_1 + \tilde\bw_2.\nonumber
\end{eqnarray}
Since $\bZ\bgamma+\bveps\sim N(0, \bP_z^{-1})$, we obtain that
$\tilde\bw_1 \sim N(0, \bH)$ with
\[
\bH= \bigl(\bZ_{\bMfrak_0^c}^T - \widetilde\bT_{12}^T
\widetilde\bT_{11}^{-1}\bZ_{\bMfrak_0}^T \bigr)
\bP_x \bP_z^{-1} \bP_x \bigl(
\bZ_{\bMfrak_0^c} - \bZ_{\bMfrak_0}\widetilde\bT_{11}^{-1}
\widetilde\bT_{12} \bigr).
\]
Note that $\bZ_{\bMfrak_0^c}$ is a block diagonal matrix, and the
$i$th block matrix has size $n_i\times(q_n-s_{2n})$. By Condition \ref
{con3}(A), it is easy to see that $\Lammax(\bZ\Gcal\bZ^T) \leq\max
_{1\leq i\leq N}\Lammax(\bZ_iG\bZ_i^T)\leq c_1s_{2n}$.
Thus, $\bP_x\bP_z^{-1}\bP_x = \bP_x(\sigma^2 \bI+ \bZ\Gcal\bZ^T)\bP_x
\leq(\sigma^2+c_1s_{2n})\bP_x^2= (\sigma^2+c_1s_{2n})\bP_x$.
Further,
it follows from $\widetilde\bT_{12} = \bT_{12}$ and
$\bZ_{\bMfrak_0}^T \bP_x\bZ_{\bMfrak_0} \leq\widetilde\bT_{11}$ that
\begin{eqnarray*}
\bH
&\leq&\bigl(\sigma^2+c_1s_{2n}\bigr)
\bigl(\bZ_{\bMfrak_0^c}^T - \widetilde\bT_{12}^T
\widetilde\bT_{11}^{-1}\bZ_{\bMfrak_0}^T \bigr)
\bP_x \bigl( \bZ_{\bMfrak_0^c} - \bZ_{\bMfrak_0}\widetilde
\bT_{11}^{-1}\widetilde\bT_{12} \bigr)
\\
& = &\bigl(\sigma^2+c_1s_{2n}\bigr) \\
&&{}\times\bigl(
\bZ_{\bMfrak_0^c}^T\bP_x \bZ_{\bMfrak
_0^c} + \widetilde
\bT_{12}^T\widetilde\bT_{11}^{-1}
\bZ_{\bMfrak
_0}^T \bP_x\bZ_{\bMfrak_0}\widetilde
\bT_{11}^{-1}\widetilde\bT_{12} - 2
\bZ_{\bMfrak_0^c}^T\bP_x\bZ_{\bMfrak_0} \widetilde
\bT_{11}^{-1}\widetilde\bT_{12} \bigr)
\\
& \leq &\bigl(\sigma^2+c_1s_{2n}\bigr) \bigl(
\bZ_{\bMfrak_0^c}^T\bP_x \bZ_{\bMfrak_0^c}-\widetilde
\bT_{12}^T\widetilde\bT_{11}^{-1}
\widetilde\bT_{12} \bigr)\leq\bigl(\sigma^2+c_1s_{2n}
\bigr)\bZ_{\bMfrak
_0^c}^T\bP_x \bZ_{\bMfrak_0^c}.
\end{eqnarray*}
Thus, the $i$th diagonal element of $\bH$ is bounded from above by the
$i$th diagonal element of $(\sigma^2 + c_1s_{2n})\bZ_{\bMfrak
_0^c}^T\bP_x \bZ_{\bMfrak_0^c}$, and is thus bounded by $\tilde
c_1s_{2n}m_n$ with $\tilde c_1$ some positive constant. Therefore by
the normality of $\tilde\bw_1$ we have
\begin{eqnarray*}
&&
P \bigl(\|\tilde\bw_1\|_\infty\geq\bigl\{2\tilde
c_1s_{2n}m_n\log\bigl(N(q_n-s_{2n})
\bigr)\bigr\}^{1/2} \bigr)
\\
&&\qquad\leq N(q_n-s_{2n})P \bigl(\bigl|N(0,\tilde
c_1s_{2n}m_n)\bigr|\geq\bigl\{2\tilde
c_1s_{2n}m_n\log\bigl(N(q_n-s_{2n})
\bigr)\bigr\}^{1/2} \bigr)
\\
&&\qquad=O \bigl(\bigl(\log\bigl(N(q_n-s_{2n})\bigr)
\bigr)^{-1/2} \bigr)=o(1).
\end{eqnarray*}
Therefore, $\|\tilde\bw_1\|_\infty= o_p (\{s_{2n}m_n\log
(N(q_n-s_{2n})\}^{1/2} )= o_p(nN^{-1/2}\lambda_n)$ and
%
\begin{equation}
\label{e040} \max_{j>s_{2n}} \Biggl\{\sum_{k=1}^N
\tilde w_{1,jk}^2 \Biggr\}^{1/2}\leq\sqrt{N}\|\tilde
\bw_1\|_\infty= o_p(n\lambda_n)=
o_p(1)np'_{\lambda_n}(0+),
\end{equation}
where $\tilde w_{1,jk}$ is the $ ((j-1)q_n + k )$th element of
$Nq_n$-vector $\tilde\bw_1$.

Now we consider $\tilde\bw_2$. Define $\tilde\bZ_j$ as the
submatrix of $\bZ$ formed by columns corresponding to the $j$th random
effect. Then, for each $j=s_{2n}+1, \ldots, q_{n}$, by Condition \ref
{con3}(A) we obtain that
\begin{eqnarray*}
\Biggl\{\sum_{k=1}^N\tilde
w_{2,jk}^2\Biggr\}^{1/2} &=& n\lambda_n
\bigl\|\tilde\bZ_j^T\bP_x \bZ_{\bMfrak_0}
\widetilde\bT_{11}^{-1}\bu\bigl(\hbgamma_{\bMfrak
_0}^*\bigr)
\bigr\|_2\\
&\leq& n\bigl\|\lambda_n\bu\bigl(\hbgamma_{\bMfrak_0}^*
\bigr)\bigr\|_2\bigl\|\tilde\bZ_j^T\bP_x
\bZ_{\bMfrak_0}\widetilde{\bT}_{11}^{-1}\bigr\|_2,
\end{eqnarray*}
where $\tilde w_{2,jk}$ is the $ ((j-1)q_n + k )$th element of
$Nq_n$-vector $\tilde\bw_2$.
Since \mbox{$\hbgamma_{\bMfrak_0}^*\in\mathcal{N}_1$}, by~(\ref{e009}),
(\ref{e056}) and the decreasing property of $p'_{\lambda_n}(\cdot)$
we have $\|\lambda_n\bu(\hbgamma_{\bMfrak_0}^*)\|_2 \leq
p'_{\lambda_n}(\sqrt{N}b_0^*/2)$. By~(\ref{e044}) in Condition \ref
{con3}(A),
\[
\max_{j\geq s_{2n}+1}\Biggl\{\sum_{k=1}^N
\tilde w_{2,jk}^2\Biggr\}^{1/2} <
np'_{\lambda_n}(0+).
\]
Combing the above result for $\tilde\bw_2$ with~(\ref{e025}) and
(\ref{e040}), we have shown that~(\ref{e014}) holds with asymptotic
probability one. This completes the proof.

\begin{supplement}
\stitle{Supplement to ``Variable selection in linear mixed effects models''}
\slink[doi]{10.1214/12-AOS1028SUPP} 
\sdatatype{.pdf}
\sfilename{aos1028\_supp.pdf}
\sdescription{We included additional simulation examples and technical
proofs omitted from the main text: simulation
Examples~A.1--A.3, and technical proofs of
Lemmas~\ref{L3}--\ref{L1} and Proposition~\ref{prop1}.}
\end{supplement}


\printaddresses

\end{document}